\documentclass{article}

\usepackage{arxiv}

\usepackage[utf8]{inputenc} 
\usepackage[T1]{fontenc}    

\usepackage{epsfig}
\usepackage{fancybox}
\usepackage{multicol}
\usepackage{color}
\usepackage{graphicx}
\usepackage{url}
\usepackage{floatrow}
\usepackage{subfig}
\usepackage{tikz}

\usepackage{booktabs}

\usepackage{amsmath}
\usepackage{amssymb}

\usepackage{amsthm}

\newtheorem{cor}{corollary}

\newtheorem{prop}{proposition}
 
\usepackage{lineno}

\makeatletter
\newcommand{\Rmnum}[1]{\expandafter\@slowromancap\romannumeral #1@}
\makeatother

\usepackage{multirow}

\usepackage{comment}

\theoremstyle{definition}

\theoremstyle{remark}

\numberwithin{equation}{section}

\title{Divergence/connection preservation scheme in the curvilinear domain with a small geometric approximation error}

\author{Sehun Chun,\\Underwood International College, Yonsei University, South Korea, \\
\texttt{sehun.chun@yonsei.ac.kr}, \\
\And 
Taejin Oh, \\Korea Institute of Atmospheric Prediction Systems, South Korea, \\
\texttt{oht@kiaps.org}}

\begin{document}
\maketitle
%

\begin{abstract}
Additional grid points are often introduced for the higher-order polynomial of a numerical solution with curvilinear elements. However, those points are likely to be located slightly outside the domain, even when the vertices of the curvilinear elements lie within the curved domain. This misallocation of grid points generates a mesh error, called \textit{geometric approximation error}. This error is smaller than the discretization error but large enough to significantly degrade a long-time integration. Moreover, this mesh error is considered to be the leading cause of conservation error. Two novel schemes are proposed to improve conservation error and/or discretization error for long-time integration caused by geometric approximation error: The first scheme retrieves the original divergence of the original domain; the second scheme reconstructs the original path of differentiation, called \textit{connection}, thus retrieving the original connection. The increased accuracies of the proposed schemes are demonstrated by the conservation error for various partial differential equations with moving frames on the sphere.
\end{abstract}



\section{Introduction}

The curvilinear mesh has been widely welcomed for various computational simulations. A curvilinear mesh consists of curved elements that are larger than the regular elements having straight edges in a highly curved boundary or domain. A larger curved element computationally allows a larger time step and a better approximation of the domain to yield better accuracy and efficiency for the numerical scheme. Thus, curvilinear elements or curved elements are widely used, particularly in the weak formulation of high-order schemes. For example, the modern Earth system model adapts weak formulations in the context of a continuous Galerkin method (also known as the  spectral element method) or discontinuous Galerkin method such as the Energy Exascale Earth System Model(E3SM) \cite{DOEE3SM}, Nonhydrostatic Unified Model of the Atmosphere (NUMA) \cite{NUMA}, Korea Meteorological Administration's KIM \cite{Choi2016}, and the Met Office's LFRic \cite{LFRic}. For models that involve a long-time integration, such as climate models, conservation properties become important. Several models such as E3SM \cite{Taylor:2010fk, Taylor:2020aa} and LFRic \cite{Cotter:2012uv, Cotter:2014vo, Melvin:2019aa} employ compatible or mimetic schemes, which use discrete operators that mimic the key mathematical properties of the divergence, gradient, and curl operators of the continuous governing equations \cite{Thuburn:2015vf} to ensure better conservation properties.

One critical issue of the curvilinear mesh in a weak formulation is that the internal grid points may not be exactly in the domain, even when the vertices and edges are exactly located in the original domain. Consider a curved line element with two vertices, $\mathbf{x}_A$ and $\mathbf{x}_B$. A mesh for this element creates the grid points $\mathbf{x}_i,~ 1 \le i \le p-1$, for a solution of polynomial order $p$. Consider the Lagrangian polynomial expansion on this curved line element as follows.

\begin{equation*}
\mathbf{x}(\xi) = \sum_{i=0}^p \mathbf{x}_i \ell(\xi), 
\end{equation*}
where the range of $\xi$ is $0 \leq \xi \leq \xi_m$ for the collocation points $\xi_m$. The $i$th order Lagrangian polynomial basis $\ell_i$ is defined as

\begin{equation*}
 \ell_i (\xi) \equiv \prod_{0 \le m \le p,~ m \neq i} \frac{ \xi - \xi_m}{ \xi_i - \xi_m}.
\end{equation*}

Suppose that the vertices are exactly in the domain, as $\mathbf{x}_0 = \mathbf{x}_A$ and $\mathbf{x}_p = \mathbf{x}_B$. The locations of the internal grid points are not analytically allocated but should be approximated in the weak formulation. Note that the internal grid points significantly affect the accuracy and stability of the weak formulation. Many schemes have been proposed to effectively place the internal grid points on a curved line in the weak formulation, such as the Galerkin projection \cite{Spencerbook} or isoparametric representation \cite{Ergatoudis1968,Spencerbook}. However, the exact allocation of the internal grid points is challenging. This misallocation of the internal grid points produces a relatively small error, called \textit{geometric approximation error} \cite{MMFNekmesh}, defined as follows.

\begin{equation}
\mbox{Geometric approximation error} \equiv \frac{1}{\sqrt{p-1}} \left [   \sum_{i=1}^{p-1} \| \mathbf{x}_i - \mathbf{x}_i^{exact} \|^2 \right ],  \label{MesherrDef}
\end{equation}
where $\mathbf{x}_i^{exact}$ is the $i$th point's exact location on the domain. Note that the mesh error is the summation of geometric approximation error and the vertex error defined as $ (1/\sqrt{2}) [  \| \mathbf{x}_0 - \mathbf{x}_0^{exact} \|$ + $\| \mathbf{x}_p - \mathbf{x}_p^{exact} \| ]$ for the curved line element.

The geometric approximation error is not significant for all the numerical schemes. The finite difference scheme with a spherical coordinate axis and corresponding metric tensor does not generate significant geometric approximation error. This is because the finite difference scheme \textit{a priori} requires the construction of all the discrete grid points. Moreover, the analytically-known metric tensor specifies the path between the discrete grid points. Thus, all the geometric information is completely determined. This contrasts with the weak formulation where the additional grid points should be located in a curved domain depending on the polynomial order of the solution. Moreover, the Piola-transformation type numerical schemes \cite{Natale, Taylor, Thuburn} are also free from geometric approximation error.

The weak formulations of high-order methods, such as spectral methods and (continuous/discontinuous) Galerkin methods, have gained broad attention for solving partial differential equations (PDEs) on curved surfaces for their distinctive advantages such as geometric flexibility, long-time stability/accuracy, and computational efficiency for a given accuracy threshold. However, it has recently been revealed that the geometric approximation error, which is unique in the weak formulation with curvilinear mesh, may degrade the overall accuracy, particularly in long-time integration \cite{MMFNekmesh}. Moreover, it has been reported that the geometric approximation error significantly affects the conservation error of mass and energy. This is because the geometric approximation error modifies the surface normal vector to change the differential operator and consequently degrades the Jacobian for accurate integration \cite{MMFNekmesh}.

One possible solution for removing the geometric approximation error is to use a high-order mesh. In the study by \cite{MMFNekmesh}, a high-order mesh, known as \textit{NekMesh} \cite{Peiro2,Peiro1,Turner}, was introduced to demonstrate the trivial geometric approximation error of a curvilinear spherical mesh. The mesh error of NekMesh converges exponentially as $p$ increases. Briefly speaking, NekMesh deforms the CAD-based low-order linear mesh of the surface to conform to the original geometry. NekMesh reallocates the interior grid points through energy-based optimization for optimal nodal sets. However, the problem is that the high-order mesh could be challenging to construct for a domain that is more complex than the sphere. Consider the diverse curved elements of a sphere having quadrilateral or hexagonal patches. A high-order mesh may not always be available with specific geometric constraints, even for the sphere. This is the case for the majority of finite-difference or finite-volume schemes on the sphere. These schemes generally show second-order convergence with grid systems that do not use the latitude-longitude grid (e.g., triangular, hexagonal, or cubed-sphere grids). This phenomenon occurs because of the lack of grid orthogonality, a tradeoff that comes to avoid pole singularity for better parallel scalability, which consequently makes high-order schemes difficult to formulate.

This paper proposes two novel schemes to significantly reduce the conservation error caused by geometric approximation error, particularly in a mesh with non-negligible geometric approximation error. The first scheme retrieves the original magnitude of the covariant divergence, and the second scheme retrieves the path of differentiation for the original covariant divergence. By covariant divergence, we mean the divergence on general curved surfaces where the vector differentiation should be considered with respect to the varying axes on the surface. 

The introduced schemes can be applied to any arbitrary curved surface. However, theoretical and test examples generally use the sphere for validation convenience. Mathematical symbols introduced in this paper is summarized in Table \ref{symbols}. Moving frames have been demonstrated to be accurate and efficient geometric tools for the numerical solution of various PDEs on general curved surfaces \cite{MMF1,MMF2,MMF3,MMF4}. Moving frames provide a convenient change of differentiation to constitute original and deformed connections. However, the error related to the different directions of moving frames on the sphere is considered negligible in this paper. This is made possible by using a sufficient number of curved elements and grid points for the domain, as described in refs. \cite{MMF1,MMF2}.

The content of this paper is organized as follows. Chapter 2 and 3 explain the changes in divergence and flux by geometric approximation error of a slightly deformed domain. In Chapter 4, divergence/connection preservation schemes are proposed to improve the conservation error of the numerical schemes for general PDEs with divergence. Chapter 5 provides explanations about how to construct LOCAL and LOCSPH moving frames on an arbitrary curved domain. Chapter 6 presents numerical tests of static divergence and curl on the sphere and corresponding partial differential equations such as conservational laws and Maxwell's equations. Chapter 7 displays several numerical tests for shallow water equations (SWEs) to show the accuracy of the proposed scheme, particularly for the conservation of mass and energy. A discussion follows in Chapter 8. 

\begin{table}[ht]
\begin{tabular}{c c c }
 \hline\noalign{\smallskip} 
 &  Original domain ($\Omega^0$) & Slightly deformed domain ($\Omega^{\varepsilon}$)  \\
  \noalign{\smallskip}\hline\noalign{\smallskip}
Surface normal vector & $\mathbf{k}^0$ & $\mathbf{k}^{\varepsilon}$  \\
Edge normal vector & $\mathbf{n}^0$ & $\mathbf{n}^{\varepsilon}$  \\
Tangent vector of the edge & $\mathbf{t}^0$ & $\mathbf{t}^{\varepsilon}$  \\
Velocity vector & $\mathbf{v}^0$ & $\mathbf{v}^{\varepsilon}$  \\
Del operator & $\nabla^0$ & $\nabla^{\varepsilon}$ \\
\hline
\end{tabular}
\caption{List of symbols}
\label{symbols}
\end{table}

\section{Divergence on a slightly deformed domain}

Consider a smoothly curved original domain $\Omega^{0}$ without any mesh error. A tessellation of the domain is obtained as $\Omega^{0} = \cup_e  \Omega_e^{0}$ where $\Omega_i^{0} \cap \Omega_j^{0} = \emptyset$ when $i \neq j$. For a velocity vector $\mathbf{v}$ defined on a curved element $\Omega_e^{0}$, the divergence theorem holds, 

\begin{equation}
\int_{\Omega_e^0} \nabla^0 \cdot \mathbf{v}^0 d \mathbf{x} = \int_{\partial \Omega_e^0} \mathbf{n}^0  \cdot \mathbf{v}^0 d s, \label{DivergenceThm0}
\end{equation}

where the superscript $0$ indicates that the corresponding quantity is defined on the original domain $\Omega^0$. Consider an element of the domain that has a geometric approximation error of the magnitude of $\varepsilon$. Let $\varepsilon$ be the angle between the surface normal vectors of the original domain and the deformed domain. Denote this element as $\Omega^{\varepsilon}_e$. We use the superscript $\varepsilon$ to represent a vector or a differential operator that lies on $\Omega^{\varepsilon}_e$. The divergence theorem still holds for the vector $\mathbf{v}^{\varepsilon}$ with the corresponding differential operator defined on $\Omega^{\varepsilon}_e$ as

\begin{equation}
\int_{\Omega^{\varepsilon}_e} \nabla^{\varepsilon} \cdot \mathbf{v}^{\varepsilon} d \mathbf{x} = \int_{\partial \Omega^{\varepsilon}_e} \mathbf{n}^{\varepsilon} \cdot \mathbf{v}^{\varepsilon} d s. \label{DivergenceThm1}
\end{equation}

On a slightly deformed surface $\Omega_e^{\varepsilon}$ with a geometric approximation error, the source and sink of the flow is conserved, respectively. The flux that enters through the boundaries of the domain is the divergence of the vector, which is often identified as the source of the flow. When the ambient space is absent, every vector is projected into the domain $\Omega_e^{\varepsilon}$. Thus, the projection of the vector is generally achieved instantaneously. 

The projection of the divergence operator $\nabla \cdot$ is critical because the path on a curved surface determines the vector differentiation between discrete grid points. To obtain accurate vector differentiation on a surface, the exact metric tensor should be provided. This means that the path between the grid points should be computed accurately. Another option is to use moving frames lying on the domain \cite{MMF3, MMF1, MMF2, MMF4}. The projection of the divergence operator $\nabla^{\varepsilon} \cdot$ is instantly obtained without computing the metric tensor of the new curved domain $\Omega_e^{\varepsilon}$. In the moving frames, the divergence operator $\nabla^{\varepsilon} \cdot$ can be conveniently constructed on the domain $\Omega_e^{\varepsilon}$. 

\begin{figure}[ht]
\centering
 \includegraphics[width=5cm]{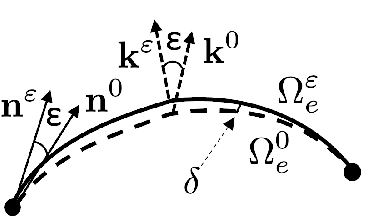}
\caption{Illustrations of the geometric components of $\Omega_e^0$ and $\Omega_e^{\varepsilon}$. }
\label {geoerror}
\end{figure}

The divergence theorems in Eq. \eqref{DivergenceThm0} and Eq. \eqref{DivergenceThm1} holds respectively. However, this does not necessarily mean that the two equalities are the same, i.e., 

\begin{equation}
\int_{\Omega^{\varepsilon}_e} \nabla^{\varepsilon} \cdot \mathbf{v}^{\varepsilon} d \mathbf{x} \neq \int_{\Omega_e^0} \nabla^0 \cdot \mathbf{v}^0  d \mathbf{x} .  \label{DivergenceComp}
\end{equation}

The divergence on a surface can be changed due to the changes in the domain. For example, consider a spherical shell of radius $r$ with the metric $ds^2 = r^2 d \theta^2 + r^2 \sin^2 \theta d \phi^2$. This implies that $x = r \sin \theta \cos \phi$, $y = r \sin \theta \sin \phi$, $z = r \cos \theta$, where $\phi$ and $\theta$ are the azimuthal and polar angle, respectively. Consider a velocity vector $\mathbf{v}^0$ on this sphere such that the velocity is expanded by the spherical coordinate axis of the sphere as 

\begin{equation*}
\mathbf{v}^0 = v_{\theta} \boldsymbol{\theta} + v_{\phi} \boldsymbol{\phi} .
\end{equation*}
where $\boldsymbol{\phi}$ and $\boldsymbol{\theta}$ represent the tangent vectors of the axis $\phi$ and $\theta$, respectively. The divergence of this vector on $\Omega$ is computed as

\begin{equation}
\nabla^0 \cdot \mathbf{v}^0 = \frac{1}{r \sin \theta} \left ( \frac{\partial v_{\phi}}{\partial \phi} + \frac{\partial }{\partial \theta} ( v_{\theta} \sin \theta )  \right ).  \label{divSphere}
\end{equation}
In the domain $\Omega^{\varepsilon}_e$, suppose that the inner grid points lie slightly outside the sphere of radius $r$. For simplicity, suppose that only the $\theta$ component changes. The velocity vector $\mathbf{v}^{\varepsilon}$ on $\Omega^{\varepsilon}_e$ is then expressed as

\begin{equation}
\mathbf{v}^{\varepsilon} =  {v}_{\theta '}   \boldsymbol{\theta '} + {v}_{\phi}    \boldsymbol{\phi} .  \label{vnewOmega}
\end{equation}
The new unit axis $\theta'$ for  $\| \boldsymbol{\theta '} \|  = 1$ is constructed as $\boldsymbol{\theta '} = \sin \varepsilon   \mathbf{r} + \cos \varepsilon  \boldsymbol{\theta}$ where $\varepsilon$ is the angle between $\theta$ and $\theta'$. Note that the velocity coefficient difference $\| {v}_{\theta} - {v}_{\theta '} \|$ corresponds to the projection error of the velocity vector and is approximately proportional to $\mathcal{O}(\varepsilon^2)$. The grid points have a small variable of radius dependent on $\theta$ such that $r' = r + \sin \varepsilon (\theta')$ where
$x = r'   \sin \theta' \cos \phi$, $y = r' \sin \theta' \sin \phi$, $z = r' \cos \theta'$. Then, the divergence of the velocity vector on $\Omega^{\varepsilon}_e$ is approximately given as 

\begin{equation}
\nabla^{\varepsilon} \cdot \mathbf{v}^{\varepsilon} = \frac{1}{ r'  \sin \theta'} \left [ \frac{\cos \varepsilon}{ r' } \frac{\partial  \varepsilon (\theta')}{\partial \theta '} v_{\theta'} \sin \theta'   + \left ( \frac{\partial v_{\phi}}{\partial \phi} + \frac{\partial }{\partial \theta'} ( v_{\theta'} \sin \theta' )  \right )  \right ] .  \label{divmodSphere}
\end{equation}

Suppose that $\varepsilon$ is sufficiently small such that $\theta \approx \theta'$. Then, for the projected velocity vector $\mathbf{v}^{\varepsilon}$, the difference between the two different divergence operator Eqs. \eqref{divSphere} and \eqref{divmodSphere} is given as 

\begin{equation}
\| \nabla^{\varepsilon} \cdot \mathbf{v}^{\varepsilon}  - \nabla^0  \cdot \mathbf{v}^0  \| \approx   \frac{v_{\theta'}  }{ r^2  }  \frac{\partial  \varepsilon (\theta')}{\partial \theta '} .  \label{divdiff}
\end{equation}

Eq. \eqref{divdiff} shows the changes in divergence, or the source term, by geometric approximation error. The change in divergence occurs due to the curvedness of the domain. Thus, it depends on the direction of the velocity vector. A slight modification of flux magnitude along a specific axis, for example, $v_{\theta'}$ in Eq. \eqref{divdiff}, dramatically changes the divergence of the flux. However, any modification of magnitude along another axis, for example, $v_{\phi}$ in Eq. \eqref{divdiff}, does not affect the divergence of the flux. The divergence also remains the same if the right-hand side of Eq. \eqref{divdiff} is zero. The change in surface area by the deformed domain causes further changes in Eq. \eqref{divdiff}. 

\section{Flux on a slightly deformed domain}              

Eq. \eqref{divdiff} indicates that the slightly deformed domain $\Omega^{\varepsilon}_e$ yields a change in the divergence of the projected vector $\mathbf{v}^{\varepsilon}$, in comparison with the divergence of the vector $\mathbf{v}^0$ in $\Omega_e^0$. Then, the inequality Eq. \eqref{DivergenceComp} corresponds to

\begin{equation}
\int_{\partial \Omega^{\varepsilon}_e} \mathbf{n}^{\varepsilon} \cdot \mathbf{v}^{\varepsilon} d s \neq \int_{\partial \Omega_e^0} \mathbf{n}^{0} \cdot \mathbf{v}^{0} d s \label{FluxComp}
\end{equation}

Note that the left and right integrand have the same magnitude. For example, consider the edge normal vector and the velocity vector lying on the sphere $\mathbf{v}$ as

\begin{equation*}
\mathbf{n}^{0} = n_{\theta} \boldsymbol{\theta} + n_{\phi} \boldsymbol{\phi}, ~~~ \mathbf{v}^{0} = v_{\theta} \boldsymbol{\theta} + v_{\phi} \boldsymbol{\phi}.
\end{equation*}

Consider a slightly deformed domain $\Omega^{\varepsilon}_e$ with an angle of $\varepsilon$, as defined previously. For simplicity, suppose that deformation occurs only along the $\theta$-axis. Similar conclusion can be drawn for the general case. Then, the unit tangent vector and the velocity vector are modified as

\begin{align*}
 \mathbf{n}^{\varepsilon} &= - n_{\theta} \sin \varepsilon \mathbf{r} + n_{\theta} \cos \varepsilon \boldsymbol{\theta} + n_{\phi} \boldsymbol{\phi}. \\
 \mathbf{v}^{\varepsilon} &= - v_{\theta} \sin \varepsilon \mathbf{r} + v_{\theta} \cos \varepsilon \boldsymbol{\theta} + v_{\phi} \boldsymbol{\phi}.
\end{align*}

Simple calculus yields the same flux as

\begin{equation}
\mathbf{n}^{0} \cdot \mathbf{v}^{0} =  n_{\theta} v_{\theta} +  n_{\phi} v_{\phi}  = \mathbf{n}^{\varepsilon} \cdot \mathbf{v}^{\varepsilon} \label{flux1}
\end{equation}
Therefore, the only other reason for the inequality in Eq. \eqref{FluxComp} is the length difference between the lines $\partial \Omega^{\varepsilon}_e$ and $\partial \Omega_e^0$. 

Let $\mathbf{k}^0$ and $\mathbf{k}^{\varepsilon}$ be the surface normal vector of $\Omega^0$ and $\Omega^{\varepsilon}$, respectively, as illustrated in Fig. \ref{geoerror}. Let $\varepsilon$ be the angle between $\mathbf{k}^0$ and $\mathbf{k}^{\varepsilon}$. The angle between the corresponding edge normal vector $\mathbf{n}^0$ and $\mathbf{n}^{\varepsilon}$ is also $\varepsilon$. For simplicity, suppose that there is no deformation along the direction orthogonal to the edge normal vector, equivalently $\mathbf{t}^0 = \mathbf{t}^{\varepsilon}$. A similar conclusion could be reached for the general case where deformation occurs both for $\mathbf{n}$ and $\mathbf{t}$. Define the edge length $h$ of the element as follows: 

\begin{equation}
h_e = \int_{\partial \Omega_e^0} ds, ~~~\mbox{or},~~~h^{\varepsilon}_e = \int_{\partial \Omega^{\varepsilon}_e} ds. \label{Defh}
\end{equation}

The boundaries $\partial \Omega^{\varepsilon}_e$ and $\partial \Omega_e^0$ become flatter as the edge length $h$ becomes smaller. Suppose that $h$ is sufficiently small. Note that the curvature  of the lines in Fig. \ref{geoerror} is exaggerated but it is common to have thousands of curved elements or more for a sphere of radius $r$. Thus, every line on each element is close to a straight line. Moreover, the angle $\varepsilon$ is sufficiently small because we assume that the geometric approximation error is relatively smaller than the discretization error. Let the mesh error, which is defined as the radius difference between $\partial \Omega^{\varepsilon}_e$ and $\partial \Omega_e^0$, be denoted by $\delta$. Then, the mesh error $\delta$ is expressed as follows under the above assumption:

\begin{equation*}
\delta \approx \frac{h}{2} \tan \varepsilon = \mathcal{O}(\varepsilon).
\end{equation*}

This means that the mesh error is the first order function of the angle $\varepsilon$. With the assumption of sufficiently small $h$, the ratio of the length between $\partial \Omega^{\varepsilon}_e$ and $\partial \Omega_e^0$ is approximated as follows:

\begin{equation}
 h^{\varepsilon}_e  \approx  \int_{\partial \Omega_e^0} \sec \varepsilon(\mathbf{x}) ds . \label{lendiff1}
\end{equation}

Particularly, if $\varepsilon$ is approximately constant within the domain, then we have 

\begin{equation}
h^{\varepsilon}_e  \approx  \frac{h}{\cos \varepsilon}. \label{lendiff2}
\end{equation}

Thus, we have the following change of flux for the slightly deformed domain $\Omega^{\varepsilon}_e$:

\begin{equation}
 \int_{\partial \Omega_e^0} \mathbf{n}^{0} \cdot \mathbf{v}^{0} d s =  \int_{\partial \Omega^{\varepsilon}_e} \cos \varepsilon \mathbf{n}^{\varepsilon} \cdot \mathbf{v}^{\varepsilon} d s \label{Fluxrelation}.
\end{equation}

Eqs. \eqref{lendiff1} and \eqref{lendiff2} explains the inequality in Eq. \eqref{FluxComp}. Therefore, from Eqs. \eqref{DivergenceComp} and \eqref{lendiff1}, we reach the following proposition.

\begin{prop}
A small geometric approximation error can generate a small change in the divergence of the flow on a curved domain. \label{prop1}
\end{prop}

This proposition suggests that the flow in a curved surface is caused by the geometric shape of the domain, particularly when the domain is curved. Consequently, the change of the shape of the domain may have a direct impact on the divergence of the flow. 

Consequently, the goal is to accurately compute the original divergence of the flux in the deformed domain with the geometric approximation error. Thus, the proposed scheme should identify the original divergence from the modified divergence of the deformed domain. This scheme computes the divergence caused by the geometric approximation error and use it to compute the original divergence of the velocity vector. For this purpose, the curl theorem is used.

\section{Divergence/connection preservation scheme}
Consider that the differential operator $\nabla^{\varepsilon}$ lies on $\Omega^{\varepsilon}_e$. For example, if differentiation is obtained on the discrete points on $\Omega^{\varepsilon}_e$, then the differential operator in the scheme is $\nabla^{\varepsilon}$. Construction of such a differential operator with moving frames will be explained in the next section. Consider a vector $\mathbf{F}$ lying on the curved element $\Omega^{\varepsilon}_e$ with geometric approximation error. Then, Stokes theorem is expressed with $\mathbf{F}$ in $\Omega^{\varepsilon}_e$ as follows: 

\begin{equation}
\iint_{\Omega^{\varepsilon}_e} ( \nabla^{\varepsilon}  \times \mathbf{F}   ) \cdot d \mathbf{S}  = \oint_{\partial \Omega^{\varepsilon}_e} \mathbf{F}  \cdot d \boldsymbol{\ell}  . \label{CurlThm}
\end{equation}

On the curved element $\Omega^{\varepsilon}_e$, two orthonormal basis vectors can be chosen such as $\mathbf{e}_{\varepsilon}^1$ and $\mathbf{e}_{\varepsilon}^2$ which lie on the tangent plane. The surface normal vector for $\Omega^{\varepsilon}_e$ can be defined as $\mathbf{k}^{\varepsilon} = \mathbf{e}^{\varepsilon}_1 \times \mathbf{e}^{\varepsilon}_2$. Then, the curl operator can be defined on $\Omega^{\varepsilon}_e$ such that the differential operator is expanded in $\mathbf{e}_{\varepsilon}^1$ and $\mathbf{e}_{\varepsilon}^2$, similar to the following gradient operator:

\begin{equation*} 
\nabla^{\varepsilon} f = ( \nabla f \cdot \mathbf{e}^1_{\varepsilon} ) \mathbf{e}^1_{\varepsilon} + ( \nabla f \cdot \mathbf{e}^2_{\varepsilon} ) \mathbf{e}^2_{\varepsilon} .
\end{equation*}

\subsection{Divergence preservation scheme}

Define the vector $\mathbf{F} $ in Eq. \eqref{CurlThm} as follows:

\begin{equation*}
 \mathbf{F} = \mathbf{F}^{0\varepsilon}  \equiv \mathbf{k}^{0}  \times \mathbf{v}^{\varepsilon}. 
 \end{equation*}

Note that the surface normal vector $ \mathbf{k}^{0}$ of the original domain $\Omega_e$ is used, instead of $\mathbf{k}^{\varepsilon}$. Substituting it into Eq. \eqref{CurlThm}, we obtain

\begin{equation}
\iint_{\Omega^{\varepsilon}_e} ( \nabla^{\varepsilon}  \cdot \mathbf{v}^{\varepsilon} )  d A -  \iint_{\Omega^{\varepsilon}_e }   \mathcal{G} (\mathbf{k}^0,\nabla^{\varepsilon},\mathbf{v}^{\varepsilon})    d A = \oint_{\partial \Omega^{\varepsilon}_e } {\mathbf{n}}^{0}   \cdot \mathbf{v}^{\varepsilon}  ds. \label{CurlThm1}
\end{equation}

The function $\mathcal{G}$ is defined as follows:

\begin{equation}
\mathcal{G} (\mathbf{k} ,\nabla, \mathbf{v} )   = \mathbf{k}   \cdot  \left [ (\mathbf{k}   \cdot \nabla )\mathbf{v}   -  (\mathbf{v}  \cdot \nabla )\mathbf{k} \right ]  .  \label{Gdef1}
\end{equation}

The function $\mathcal{G} $ measures the variation in $\mathbf{v} $ with respect to $\mathbf{k} $ along the paths of $\nabla $, known as the \textit{connection}. Note that $\mathcal{G}(\mathbf{k}^0,\nabla^{\varepsilon},\mathbf{v}^{\varepsilon}) $ is always zero if $\mathbf{k}^{\varepsilon}$ is used instead of $\mathbf{k}^{0}$ in Eq. \eqref{Gdef1}. Generally speaking, if $\mathbf{k}$, $\mathbf{v}$, and $\nabla$ lie on the same domain, then $\mathcal{G}$ is always zero and, consequently, the conservation laws hold in Eq. \eqref{CurlThm1}. The definition of $\mathbf{F}^{0\varepsilon}$ yields the change of the flux in Eq. \eqref{CurlThm1}, compared with the flux with $\mathbf{F}^{\varepsilon \varepsilon} = \mathbf{k}^{\varepsilon} \times \mathbf{v}^{\varepsilon}$. Consider that $\mathbf{n}^{0}   \cdot \mathbf{v}^{\varepsilon}  =  \cos \varepsilon ( \mathbf{n}^{\varepsilon}  \cdot \mathbf{v}^{\varepsilon} )$.

From Eq. \eqref{Fluxrelation}, we observe that $\mathcal{G}(\mathbf{k}^0,\nabla^{\varepsilon},\mathbf{v}^{\varepsilon})$ is the divergence difference of the velocity vector between the original domain and the deformed domain. Thus, we obtain

\begin{equation}
\iint_{\Omega^{\varepsilon}_e} ( \nabla^{0}  \cdot \mathbf{v}^{0} )  d A = \iint_{\Omega^{\varepsilon}_e} ( \nabla^{\varepsilon}  \cdot \mathbf{v}^{\varepsilon} )  d A -  \iint_{\Omega^{\varepsilon}_e }  \mathcal{G}(\mathbf{k}^0,\nabla^{\varepsilon},\mathbf{v}^{\varepsilon})  d A  , \label{Restorescheme1}
\end{equation}

Thus, we call Eq. \eqref{Restorescheme1} the \textit{divergence preservation scheme}. Consequently, we have the following corollary. 

\begin{cor}
The original divergence of the velocity vector $\mathbf{v}^0$ can be approximately obtained by adding $- \iint_{\Omega^{\varepsilon}_e }  \mathcal{G}(\mathbf{k}^0,\nabla^{\varepsilon},\mathbf{v}^{\varepsilon}) d A $, defined in Eq. \eqref{Gdef1}, from the modified divergence $\iint_{\Omega^{\varepsilon}_e} ( \nabla^{\varepsilon}  \cdot \mathbf{v}^{\varepsilon} )  d A$  on the deformed domain to the geometric approximation error. 
\label{corol1}
\end{cor}

\subsection{Connection preservation scheme}

As an alternative scheme to approximately compute the original divergence, the differential operator $\nabla^0$ of the original domain can be used with the curl theorem. Define the vector $\mathbf{F}^{\varepsilon \varepsilon} $ in Eq. \eqref{CurlThm} as follows

\begin{equation*}
 \mathbf{F}^{\varepsilon \varepsilon}  = \mathbf{k}^{\varepsilon}  \times \mathbf{v}^{\varepsilon}.    \label{F2}
 \end{equation*}  

Substituting this vector into the curl theorem in Eq. \eqref{CurlThm}, we obtain

\begin{equation}
\iint_{\Omega^{\varepsilon}_e} ( \nabla^0  \cdot \mathbf{v}^{\varepsilon} )  d A +  \iint_{\Omega^{\varepsilon}_e }   \mathcal{G}(\mathbf{k}^{\varepsilon},-\nabla^{0},\mathbf{v}^{\varepsilon})    d A = \oint_{\partial \Omega^{\varepsilon}_e } {\mathbf{n}}^{\varepsilon}  \cdot \mathbf{v}^{\varepsilon}  ds. \label{CurlThm3}
\end{equation}

The differential operator $\nabla^0$ is constructed from the original domain $\Omega_e^0$. Thus, $ \mathcal{G}(\mathbf{k}^{\varepsilon},\nabla^{0},\mathbf{v}^{\varepsilon}) $ is not always zero and can be non trivial. Eq. \eqref{CurlThm3} makes it possible to compute the divergence of $\nabla^0 \cdot \mathbf{v}^{\varepsilon}$ which does not hold conservation properties, i.e.,

\begin{equation}
\iint_{\Omega^{\varepsilon}_e} ( \nabla^0  \cdot \mathbf{v}^{\varepsilon} )  d A  = \iint_{\Omega^{\varepsilon}_e} ( \nabla^{\varepsilon}  \cdot \mathbf{v}^{\varepsilon} )  d A -  \iint_{\Omega^{\varepsilon}_e }   \mathcal{G}(\mathbf{k}^{\varepsilon}, -\nabla^{0},\mathbf{v}^{\varepsilon})    d A .  \label{Restorescheme2}
\end{equation}

The negative gradient is used to construct a divergence form similar to that in Eq. \eqref{Restorescheme1}. In summary, we have the following corollary. 

\begin{cor}
The original divergence $\iint_{\Omega^{0}_e} ( \nabla^{0}  \cdot \mathbf{v}^{\varepsilon} )  d A $ of the vector $\mathbf{v}^{\varepsilon}$ can be approximately obtained by adding $\iint_{\Omega^{\varepsilon}_e }   \mathcal{G}(\mathbf{k}^{\varepsilon},\nabla^{0},\mathbf{v}^{\varepsilon})    d A $ to the modified divergence $\iint_{\Omega^{\varepsilon}_e} ( \nabla^{0}  \cdot \mathbf{v}^{\varepsilon} )  d A$  on the deformed domain with geometric approximation error. 
\label{corol2}
\end{cor}

\subsection{Computation of $\mathcal{G}$}
The modifications of the divergence shown in Eqs. \eqref{Restorescheme1} and \eqref{Restorescheme2} can be written as the approximation of the original divergence in the deformed domain $\Omega^{\varepsilon}_0$ as follows.

\begin{equation*}
\mbox{Original divergence}= \iint_{\Omega^{\varepsilon}_e} ( \nabla^{\varepsilon}  \cdot \mathbf{v}^{\varepsilon} )  d A -  \iint_{\Omega^{\varepsilon}_e }   \mathcal{G}   d A ,
\end{equation*}
where $\mathcal{G}$ is defined in Eq. \eqref{Gdef1}. The computation of the spurious divergence $\mathcal{G}$ requires a novel numerical scheme to compute material derivative, or convective operator, which is mathematically expressed as $(\mathbf{A} \cdot \nabla) \mathbf{B}$. Consider a two-dimensional orthogonal curved domain with the constant magnitude of tangent vectors. Then, the material divergence on this domain is given as follows \cite{Weisstein}.

\begin{equation}
  [ ( \mathbf{A} \cdot \nabla ) \mathbf{B} ]_i  =   \sum_{k=1}^{2}  \frac{A_k}{\sqrt{g_{kk}}} \frac{\partial B_i}{\partial \xi_k}  ,\label{matder}
\end{equation}
where $g_{kk}$ is the metric tensor and $\xi_i$ is the orthogonal curved axis. The subscript $i$ implies that the corresponding quantity is the vector expansion along the $i$th axis. Construct moving frames such that $g_{kk} = 1, ~ 1 \le k \le 3$ everywhere. Expand the vectors $\mathbf{A}$ and $\mathbf{B}$ on the frames $\mathbf{e}^i,~ 1 \le i \le 3$ such that

\begin{equation*}
\mathbf{A} = A_1 \mathbf{e}^1 + A_2 \mathbf{e}^2 + A_3 \mathbf{e}^3, ~~~~ \mathbf{B} = B_1 \mathbf{e}^1 + B_2 \mathbf{e}^2 + B_3 \mathbf{e}^3. 
\end{equation*}
Then, the numerical scheme with moving frames to compute the material derivative in Eq. \eqref{matder} is given as follows.
\begin{equation}
   ( \mathbf{A} \cdot \nabla ) \mathbf{B}  =  \sum_{i=1}^3 \left [ A_1 (\nabla B_i \cdot \mathbf{e}^1 ) +  A_2 (\nabla B_i \cdot \mathbf{e}^2) \right ] \mathbf{e}^i . \label{matdermf}
\end{equation}
Then, the spurious divergence $\mathcal{G}$ is computed as follows:
\begin{equation}
\mathcal{G} =  \sum_{i=1}^3 k_i \left [ k_1 (\nabla v_i \cdot \mathbf{e}^1 ) +  k_2 (\nabla v_i \cdot \mathbf{e}^2) - v_1 (\nabla k_i \cdot \mathbf{e}^1 ) -  v_2 (\nabla k_i \cdot \mathbf{e}^2) \right ] . \label{Gmf}
\end{equation}

Note that if $\mathbf{k}$ is orthogonal to $\mathbf{e}^1$ and $\mathbf{e}^2$, then the spurious divergence $\mathcal{G} =0$. However, $\mathcal{G}$ can be significant. For a nontrivial $ \mathcal{G}(\mathbf{k}^0,\nabla^{\varepsilon},\mathbf{v}^{\varepsilon})$, $\mathbf{k}^0$ should not be orthogonal to $\mathbf{e}^1_{\varepsilon}$ and $\mathbf{e}^2_{\varepsilon}$. For a nontrivial $ \mathcal{G}(\mathbf{k}^{\varepsilon},\nabla^{0},\mathbf{v}^{\varepsilon}) $, $\mathbf{k}^{\varepsilon}$ should not be orthogonal to $\mathbf{e}^1_{0}$ and $\mathbf{e}^2_{0}$.

\section{Constructing $\nabla^0$ or $\nabla^{\varepsilon}$ with moving frames}

To construct two types of differential operators, i.e., $\nabla^0$ or $\nabla^{\varepsilon}$, the local orthonormal basis sets, called \textit{moving frames}, are used. Suppose $\mathbf{e}^i, ~ 1 \le i \le 3,$ is an orthonormal basis set defined on every grid point such that

\begin{equation*}
\mathbf{e}^i \cdot \mathbf{e}^j = \delta^i_j,~~~~ \| \mathbf{e}^i \| = 1 ,
\end{equation*}
where $\delta^i_j$ is the Kronecker delta. Let the basis vector, also called as \textit{frame}, $\mathbf{e}^1$ and $\mathbf{e}^2$ lie on the tangent plane of the domain. Consequently, the third frame $\mathbf{e}^3$ is aligned along the surface normal vector. The gradient operator expanded on this moving frame is defined as follows

\begin{equation*}
\nabla f = f_1 \mathbf{e}^1 +f_2 \mathbf{e}^2, 
\end{equation*}
where $ f_i = \nabla f \cdot \mathbf{e}^i$ for $i =1,~2$. Similar definitions are applied to other differential operators. Moving frames are directly related to the direction of differentiation. By defining the domain where $\mathbf{e}^1$ and $\mathbf{e}^2$ lie, $\nabla^{\varepsilon}$ and $\nabla^{0}$ can be fully represented by moving frames.

\subsection{Constructing moving frames for $\nabla^{\varepsilon}$}
The differential operator $\nabla^{\varepsilon}$ is based on the orthonormal bases that lie on the tangent plane of the surface, which are called LOCAL moving frames \cite{MMFCov}. Consider the element-wise mapping from the standard element to the curved element of the domain. Then, construction of these bases can be easily achieved by using a differentiation of the map as follows. For a more detailed description of construction, please refer to \cite{MMF1, MMF2}.

\begin{figure}[ht]
\centering
 \includegraphics[width=8cm]{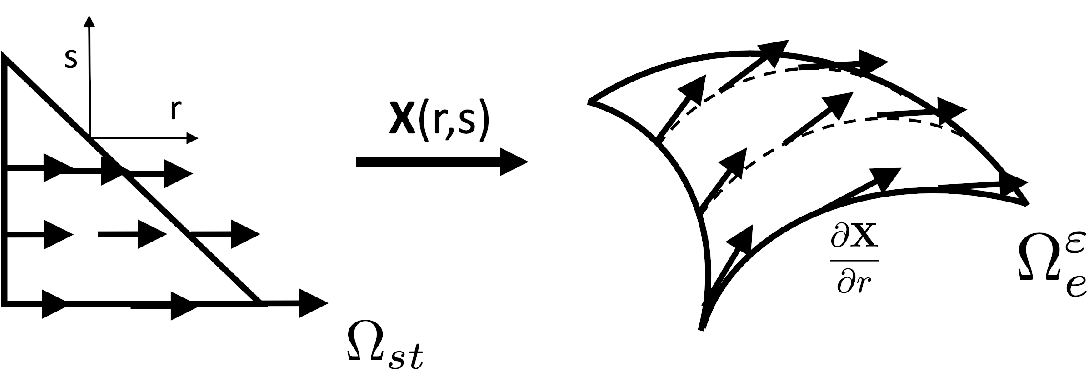}
\caption{Construction of LOCAL moving frames using the map from the standard element to a curved element.}
\label {cmap}
\end{figure}

Denote $\Omega_{st}$ as a standard quadrilateral element having unit length. Suppose $r$ and $s$ are the two independent axes of $\Omega_{st}$ as shown in Fig. \ref{cmap}. Consider a mapping $\mathbf{X}(r, s)$ from $\Omega_{st}$ to a curved element $\Omega^{\varepsilon}_e$. The mapping $\mathbf{X}(r, s)$, represented as a polynomial of two variables $r$ and $s$, displays the coordinate of the curved element $\Omega^{\varepsilon}_e$ where the Jacobian is not constant. In this mapping, the grid points and the tangent vectors of $\Omega_{st}$ correspond to the grid points and the tangent vectors of $\Omega^{\varepsilon}_e$, respectively. The tangent vector of the domain $\Omega_{st}$ is denoted by $\partial / \partial r$ or $\partial / \partial s$. Then, the mapping $\mathbf{X}(r, s)$ also maps these tangent vectors onto $\Omega^{\varepsilon}_e$ as follows: 

\begin{equation*}
\mathbf{X} \left ( \frac{\partial}{\partial r} \right ) = \frac{ \partial \mathbf{X}}{\partial r}  ~~~ \in ~~ \Omega^{\varepsilon}_e ,
\end{equation*}
where the equality holds by the law of vector mapping. The above equation implies that the tangent vector of the curved element is obtained by differentiating the map $\mathbf{c}(r, s)$ with respect to the axis. The derived two tangent vectors are not necessarily orthonormal. Thus, the orthonormalization of the two vectors yields the local orthonormal basis vectors $\mathbf{e}^1$ and $\mathbf{e}^2$. The differentiability of these bases is guaranteed by sufficient smoothness of the map $\mathbf{X}(r,s)$. The third basis $\mathbf{e}^3$ is obtained by the curl of the two basis vectors. The obtained LOCAL moving frames are generally discontinuous across the interfaces of the elements. An additional process could be applied, though it is not necessary, to construct a more Euclidean-like moving frames by comparing the divergence of $\mathbf{e}^1$ along each edge direction \cite{MMFCov}. The LOCAL moving frames coincide with the local shape of the surface. Thus, any differential operator using LOCAL moving frames is $\nabla^{\varepsilon}$.

\subsection{Constructing moving frames for $\nabla^{0}$}
Deriving the differential operator $\nabla^{0}$ defined on the original domain is also conveniently obtained by using moving frames. The only adjustment needed from LOCAL moving frames is to align $\mathbf{e}^3$ along the surface normal vector of the original domain $\Omega^0$ and reorthornomalize the frames. The first step is to replace the third frame $\mathbf{e}^3$ with the unit surface normal vector of $\Omega^0$, or $\mathbf{k} / \| \mathbf{k} \|$. The second step is to project $\mathbf{e}^1$ and $\mathbf{e}^2$ onto the new tangent space to form a new orthonormal set $\mathbf{e}^i~1, \le i \le 3$. Another orthonormalization of the basis vectors creates a set of LOCAL moving frames without geometric approximation error. When the surface normal vector is the same as that of a sphere with unit length, moving frames are called LOCSPH \textit{(short for LOCAL SPHERE)} moving frames in this paper. The third moving frame $\mathbf{e}^3$ of LOCSPH is the same as the surface normal vector of the sphere. However, major difference is that LOCSPH is discontinuous across the interfaces of elements, in contrast to the continuous spherical coordinate system. Consequently, no geometric singularity exists for LOCSPH. A quick summary of the three different moving frames is given in Table. \ref{Comparison}

\begin{table}[ht]
\begin{tabular}{c c c c}
 \hline\noalign{\smallskip} 
 & Spherical coordinate axis  & LOCSPH & LOCAL  \\
  \noalign{\smallskip}\hline\noalign{\smallskip}
Surface normal vector & $\mathbf{r}^{(*)}$ & $\mathbf{r}^{(*)}$ & Analytically Unknown \\
Continuity$^{(**)}$ & Yes & No  & No \\
\hline
\end{tabular}
\caption{Comparison of the spherical coordinate system, LOCSPH, and LOCAL moving frames, on the sphere. (*) $\mathbf{r}$ is the radial vector of the spherical coordinate system. (**) Continuity concerning the continuity of frames across the interfaces of elements. }
\label{Comparison}
\end{table}

\section{Numerical tests} 
For the test problems, the following spherical mesh known as \textit{ProjMesh} \cite{risser} is used. The cube edges of unit length are equidistantly dissected according to a user-defined edge length. As ProjMesh projects the mesh of the cube onto a perfect sphere, the vertices of the derived edges are exactly on the sphere. Vertex and node distributions are uniform. However, the edges of triangular elements are inaccurately represented. Consequently, the inaccurate edges yields inaccurate locations for internal grid points by isoparametric representation \cite{Spencerbook,Ergatoudis1968}.

Consider the definition of mesh error and surface normal vector error in $L_2$ as follows. Let $r^{\textrm{0}}$ and $ \mathbf{k}^{\textrm{0}}$ be the radius and surface normal vector of the original domain, respectively. For all the indices $i$ of the total $N$ grid points in the slightly deformed domain $\Omega^{\varepsilon}$, we have
\begin{linenomath}
\begin{equation}
\mathrm{Mesh~error} =  \frac{1}{\sqrt{N}}  \sqrt{ \sum_{i=1}^{N} {( r^0 - {r}^{\varepsilon}_i )^2 } },  \label{defmesherr}
\end{equation}
and
\begin{equation}
\mbox{Surface normal vector error} = \frac{1}{\sqrt{N}} \sqrt{ \sum_{i=1}^{N} { \| \mathbf{k}^0 - \mathbf{k}_{i}^{\varepsilon} \|^2 } } . \label{defSNerr}
\end{equation}
\end{linenomath}

\begin{figure}[ht]
\centering
\subfloat[Mesh error] {\label{SWERossbyInit}\includegraphics[width=5cm]{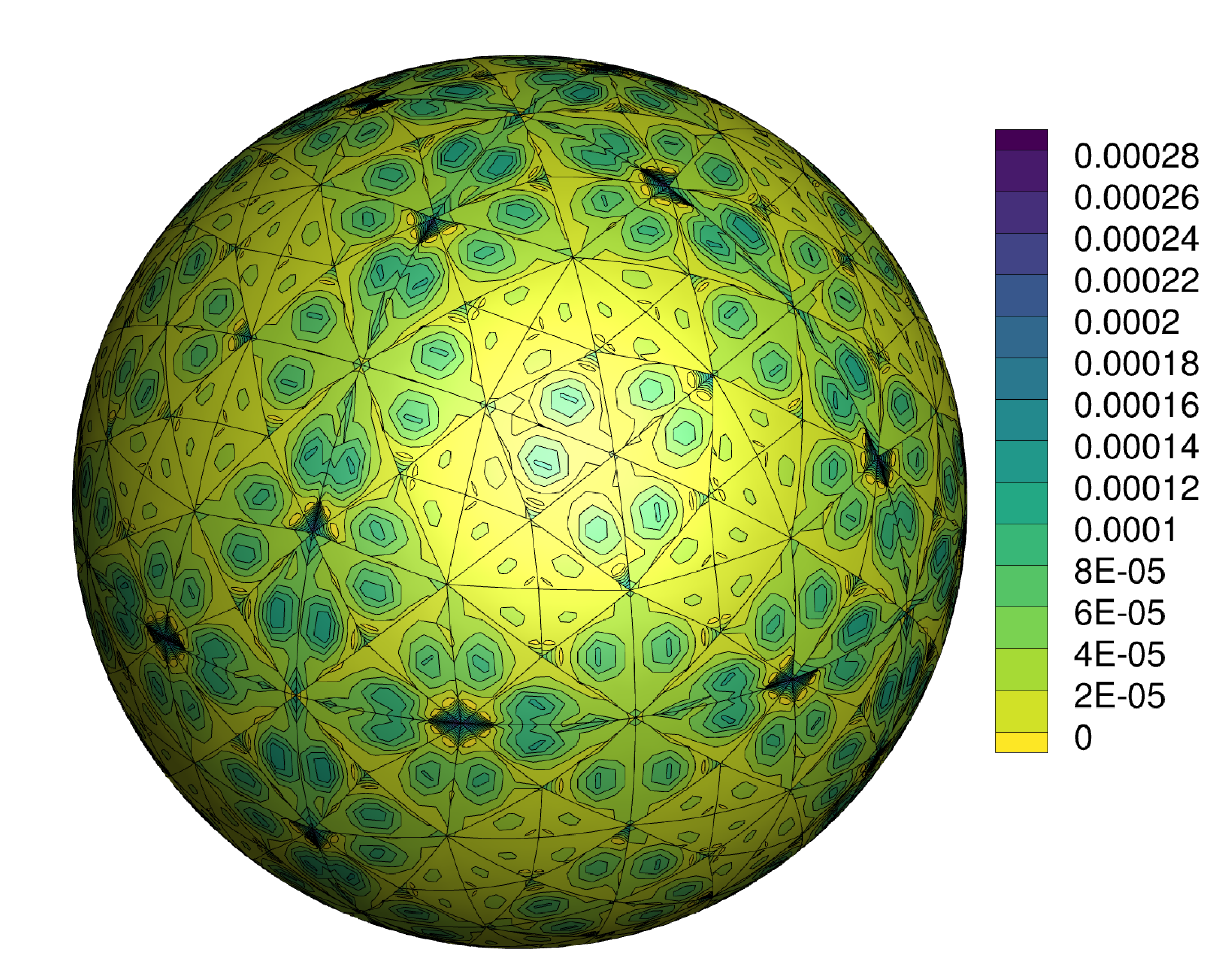} } \hspace{1cm}
\subfloat[Surface normal vector error] {\label{SWERossbyError}\includegraphics[width=5cm]{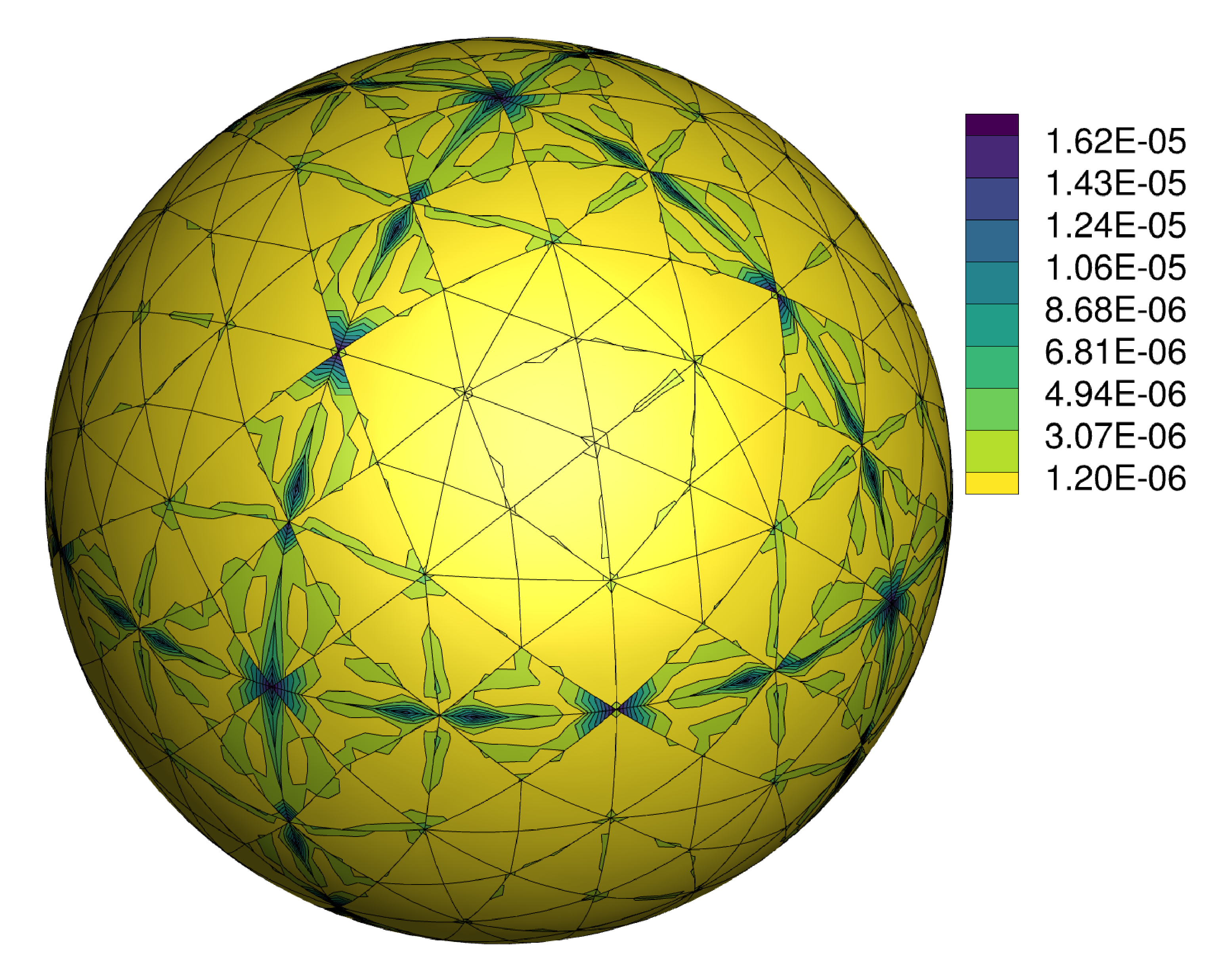} }
\caption{Distribution of the mesh error and surface normal vector error for the mesh with geometric approximation error. }
\label {Fig:Mesherr}
\end{figure}

The distribution of mesh error and the surface normal vector are displayed in Fig. \ref{Fig:Mesherr}. The maximum edge length is $h=0.4$ with 480 curved elements. Each edge has the maximum polynomial order of three. Table \ref{Tab:Mesherr} shows the convergence of mesh error as $p$. Observe that there is little convergence for the $p\ge 3$ of mesh error due to geometric approximation error. 

\begin{table}[ht]
\begin{tabular}{c c c c c c}
 \hline\noalign{\smallskip} 
 p  & 2 & 3  & 4  &  5  &  6 \\
  \noalign{\smallskip}\hline\noalign{\smallskip}
$L_2$ & 2.0497e-05 &  3.315e-05   &  3.093e-05  &  3.060e-05   & 3.045e-05      \\
$L_{inf}$  &7.565e-05 &  1.165e-04  &  1.113e-04 &  1.262e-04  &  1.098e-04     \\
\hline
\end{tabular}
\caption{Mesh error of ProjMesh with $h=0.4$. Number of curved triangular elements = 480.}
\label{Tab:Mesherr}
\end{table}

\subsection{Divergence test} 

Consider the following static divergence test on a sphere of unit length. Let $(x,y,z)=(\sin \theta \cos \phi,\sin \theta \sin \phi,\cos \theta)$. Expand the velocity vector along spherical coordinate axes, as $\mathbf{v} = v_{\theta} \boldsymbol{\theta} + v_{\phi} \boldsymbol{\phi}$. For the velocity vector, consider the Rossby-Haurwitz flow defined as
\begin{align*}
v_{\phi} &=  \omega \sin \theta + K \sin^3 \theta ( 4 \cos^2 \theta - \sin^2 \theta ) \cos 4 \phi, \\
v_{\theta} &=  -4 K  \sin^3 \theta \cos \theta \sin 4 \phi , 
\end{align*}
where $\omega = K = 7.848 \times 10^{-6}~ s^{-1}$. The divergence of the Rossby-Haurwitz flow on the sphere is zero such as

\begin{equation*}
\int_{\Omega} \nabla \cdot \mathbf{v}^i = 0,  
\end{equation*}
By using Eq. \eqref{CurlThm1}, a discontinuous Galerkin scheme for the above divergence with respect to a test function $\varphi$ is expressed as

\begin{equation}
\int_{\Omega} ( \nabla^{\varepsilon} \cdot \mathbf{v}^{\varepsilon} ) \varphi dA = - \int_{\Omega} \nabla^{\varepsilon} \varphi \cdot \mathbf{v}^{\varepsilon} d A + \int_{\partial \Omega} ( \mathbf{n}^{\varepsilon} \cdot \mathbf{v}^* ) \varphi d \ell + \int_{\Omega} \mathcal{G} \varphi dA,  \label{DGDivscheme}
\end{equation}
where $*$ represents the corresponding quantity that is the numerical flux at the interfaces. Upwind flux is chosen for accurate integration. 

\begin{figure}[ht]
\centering
\includegraphics[width=12cm]{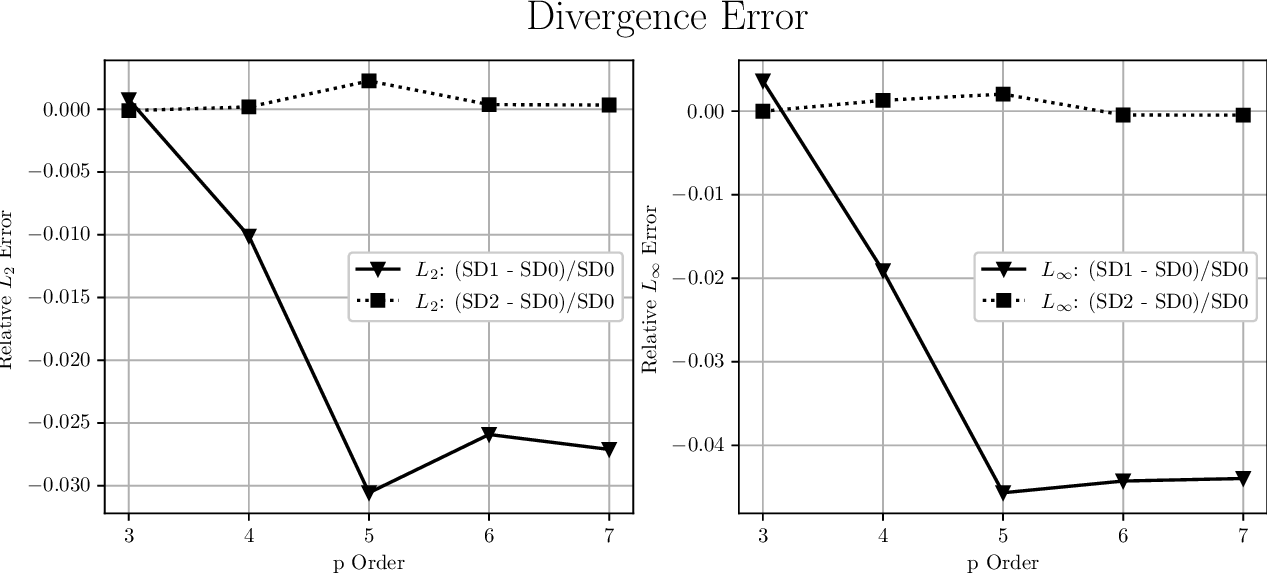}
\caption{Relative $L_2$ and $L_{\infty}$ error in static divergence test problem: SD0: LOCAL moving frames with $\mathcal{G}=0$. Divergence conservation scheme (SD1): Use Eq. \eqref{Restorescheme1}, Connection preservation scheme (SD2): Use Eq.\eqref{Restorescheme2}.}
\label {Divresult}
\end{figure}

Fig. \ref{Divresult} demonstrates enhanced accuracy of the divergence preservation scheme (SD1) or both in $L_2$ and $L_{\infty}$ error, compared with the scheme having LOCAL moving frames. Th connection preservation scheme (SD2) does not show any corresponding advantage for $L_2$ and $L_{\infty}$ error. There is no improvement in the accuracy of SD2 because the weak formulation of the divergence in Eq. \eqref{DGDivscheme} does not require any covariant differentiation. Only divergence compensation increases the accuracy of the numerical scheme noticeably. 

\subsection{Time-dependent advection}
For a velocity vector $\mathbf{v}$ on the sphere ($S^2$), the following conservation law is considered.

\begin{equation}
\frac{\partial u}{\partial t} + \nabla \cdot ( u \mathbf{v}  ) = 0. ~~~ \label{CL}
\end{equation}
Then, the weak formulation of the above conservation law is obtained as 

\begin{equation}
\int_{S^2} \frac{du}{dt} \varphi  - \int_{S^2} \nabla \varphi \cdot   \mathbf{v} u d A + \int_{\partial S^2} ( \mathbf{n} \cdot   \mathbf{v}^* ) u \varphi d \ell - \int_{S^2} \mathcal{G}  \varphi dA = 0.  \label{DGCLscheme}
\end{equation}

A smooth cosine bell is advected around the sphere to test the above conservation law. Refer to \cite{MMF1, Janusz2011, Williamson1992} for a detailed description of the test  The test problem is to passively advect the cosine bell one revolution around the sphere. Time is marched using a step size of $10^{-4}$ by an explicit fourth-order Runge-Kutta scheme. The cosine bell is advected at an angle of $\pi/4$ with respect to the axes of the poles. 

\begin{figure}[ht]
\centering
\includegraphics[width=12cm]{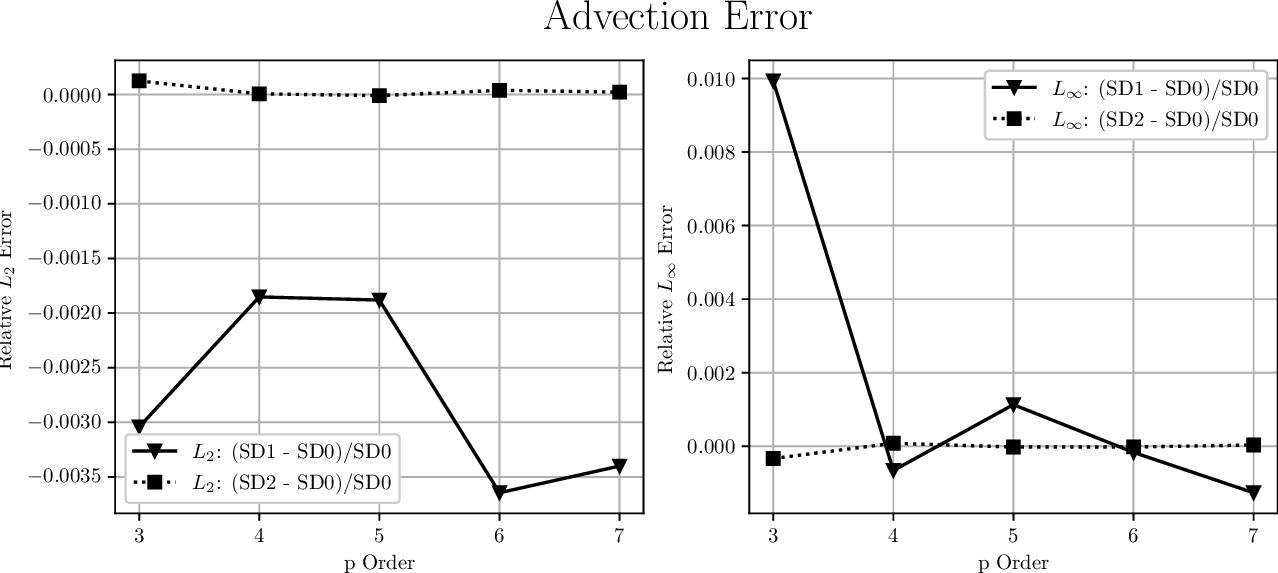}
\caption{Relative $L_2$ and $L_{\infty}$ error in time-dependent advection test problem: SD0: LOCAL moving frames with $\mathcal{G}=0$. Divergence conservation scheme (SD1) uses Eq. \eqref{Restorescheme1}. Connection preservation scheme (SD2) uses Eq.\eqref{Restorescheme2}.}
\label {Advresult}
\end{figure}

Fig. \ref{Advresult} presents an error convergence similar to that of the divergence test in Fig. \ref{Divresult}. For $L_2$ error, SD1 is significantly more accurate than the scheme having LOCAL moving frames or the SD2 scheme. For $L_{\infty}$ error, SD1 fluctuates for a low order $p$, but converges to show better accuracy for $p\ge6$.

\subsection{Curl Test} 
The cur test problem for the curl with respect to the surface normal direction $( \nabla \times \mathbf{v} ) \cdot \mathbf{k}$ is similar to the divergence test. For the test, a flipped Rossby Haurwitz flow between the $\theta$-axis and $\phi$ axis is used, and it is defined as

\begin{align*}
v_{\theta} &=  \omega \sin \theta + K \sin^3 \theta ( 4 \cos^2 \theta - \sin^2 \theta ) \cos 4 \phi, \\
v_{\phi} &=  4 K  \sin^3 \theta \cos \theta \sin 4 \phi , 
\end{align*}
where $\omega$ and $K$ are the same as that of the divergence test. Then, the both flows are curl-free to have the following equality

\begin{equation}
\int_{\Omega} ( \nabla \times \mathbf{v} ) \cdot \mathbf{k} \varphi d A = 0. \label{curlcdotk}
\end{equation}
The geometric effects on $( \nabla \times \mathbf{v} ) \cdot \mathbf{k}$ are more complicated than the divergence of a vector because the magnitude of the curl is measured only along the surface normal vector. For example, the LOCAL and LOCSPH have different surface normal directions when their tangent planes are different. A discontinuous Galerkin scheme to compute Eq. \eqref{curlcdotk} is given as

\begin{align}
\lefteqn{ \int_{\Omega} ( \nabla \times \mathbf{v} ) \cdot \mathbf{e}^3 \varphi dA = - \int_{\Omega} \nabla \varphi \cdot \tilde{\mathbf{v}}^c d A + \int_{\partial \Omega} ( \mathbf{n} \cdot \tilde{\mathbf{v}}^{c*} ) \varphi d \ell } \hspace{6cm} \nonumber \\
& + \int v^i \mathbf{e}^i \cdot ( \nabla \times \mathbf{e}^3 ) \varphi d x - \int_{\Omega} \mathcal{G}   \varphi dA,  \label{DGCurlscheme}
\end{align} 
where $\tilde{\mathbf{v}}^c \equiv - (\mathbf{v} \cdot \mathbf{e}^1) \mathbf{e}^2 +  (\mathbf{v} \cdot \mathbf{e}^2)  \mathbf{e}^1$. 

\begin{figure}[ht]
\centering
\includegraphics[width=12cm]{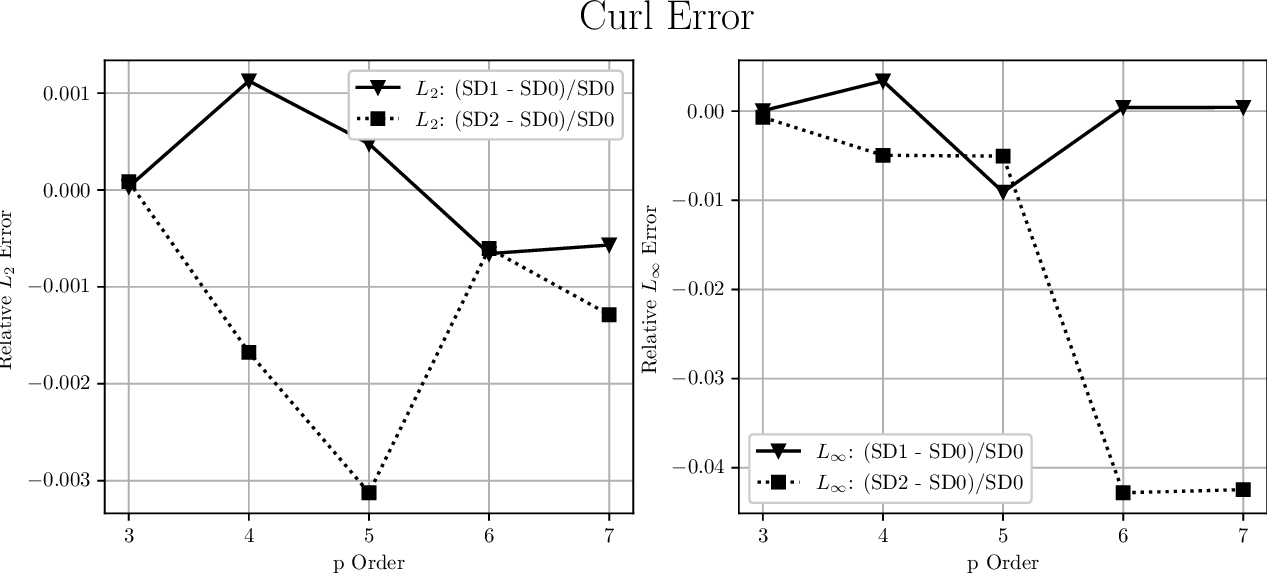}
\caption{Relative $L_2$ and $L_{\infty}$ error in static curl test problem: SD0: LOCAL moving frames with $\mathcal{G}=0$. Divergence conservation scheme (SD1) uses Eq. \eqref{Restorescheme1}. Connection preservation scheme (SD2) uses Eq.\eqref{Restorescheme2}.}
\label {Curlresult}
\end{figure}

Fig. \ref{Curlresult} demonstrates superior accuracy of the connection preservation scheme for the curl operator than that of the conservation preservation scheme. For $L_2$ error, SD2 is significantly more accurate than the scheme having LOCAL moving frames or the SD1 scheme. However, the difference between the two schemes becomes smaller as $p$ increases. For $L_{\infty}$ error, SD1 shows similar accuracy up to $p=5$. For $p \ge 6$, the SD2 scheme shows smaller error than the SD1 scheme. Overall, the SD1 scheme shows better accuracy than the scheme having LOCAL moving frames with no $\mathcal{G}$, particularly as $p$ increases. However, the best accuracy is obtained by the SD2 scheme. This is due to the covariant differentiation term, i.e., $\nabla \times \mathbf{e}^3$, which significantly contributes to the overall error of the curl operator in Eq. \eqref{DGCurlscheme}. Thus, error caused by wrong connection is greater than that caused by wrong divergence in the overall error. Consequently, retrieving the original connection improves the accuracy of the scheme. This phenomenon is more apparent at a lower $p$ because the connection error diminishes as $p$ increases.

\subsection{Time-dependent Maxwell's equations}
For time-dependent PDEs with the curl operator, the following Maxwell's equations are chosen:
\begin{equation}
\hat{\varepsilon} \frac{\partial \mathbf{E}}{\partial t} = \nabla \times \mathbf{H}, \hspace{0.5cm} \hat{\mu} \frac{\partial \mathbf{H}}{\partial t} = - \nabla \times \mathbf{E},  \label{Maxwells}
\end{equation}
where $\mathbf{E}$ and $\mathbf{H}$ are the electric field and $h$-field, respectively. $\hat{\varepsilon}$ and $\hat{\mu}$ are the permittivity and permeability tensors, respectively. Consider transverse magnetic (TM) mode such that $\mathbf{H} = H^1 \mathbf{e}^1 + H^2 \mathbf{e}^2$ and $\mathbf{E} = E^3 \mathbf{e}^3$. The transverse electric (TE) mode would yield similar results. Then, the weak formulation of Maxwell's equations in Eq. \eqref{Maxwells} in the context of discontinuous Galerkin methods is given as follows: For $i=1,2$, we have

\begin{align}
\lefteqn{ \int \mu^i \frac{\partial H^i}{\partial t} \varphi dx = \int \nabla \varphi \cdot \mathbf{E}^*_i dx - \int \mathbf{E} \cdot ( \nabla \times \mathbf{e}^i ) \varphi d x } \hspace{6cm} \nonumber \\
& + \int_{\partial \Omega} \mathbf{n} \cdot \widetilde{\mathbf{E}}^*_i \varphi d s - \int \mathcal{G}   \varphi dx , \label{Hfield}
\end{align}
where $\mathbf{E} = E^3 \mathbf{e}^3$ and $\mathbf{E}^*_i = {E}^3 \mathbf{e}^{3i} = E^3 ( \mathbf{e}^3 \times \mathbf{e}^i )$. The tilde sign implies that the corresponding quantity is an approximated value at the boundary, called numerical flux. The permeability tensor is also decomposed into two frames $\mathbf{e}^1$ and $\mathbf{e}^2$ as $\hat{\mu} = \mu^1 \mathbf{e}^1 + \mu^2 \mathbf{e}^2$. The electric field is updated as follows.

\begin{align}
\lefteqn{\int \varepsilon^3 \frac{\partial E^3}{\partial t} \varphi dx =  \sum_{i=1}^2 \left [ - \int \nabla \varphi \cdot \mathbf{H}^*_i dx + \int \mathbf{H}^i \cdot ( \nabla \times \mathbf{e}^3) \varphi d x \right . } \hspace{6cm} \nonumber \\
&+  \left .  \int_{\partial \Omega} \mathbf{n} \cdot \widetilde{\mathbf{H}}^*_i \varphi d s + \int \mathcal{G} (\mathbf{H}^*_i) \varphi dx \right ] , \label{Efield}
\end{align}
where $\mathbf{H}^i = H^i \mathbf{e}^i$ and $\mathbf{H}^*_i = {H}^i \mathbf{e}^{i3} =  H^i ( \mathbf{e}^i \times \mathbf{e}^3 )$.

\begin{figure}[ht]
\centering
\includegraphics[width=8cm]{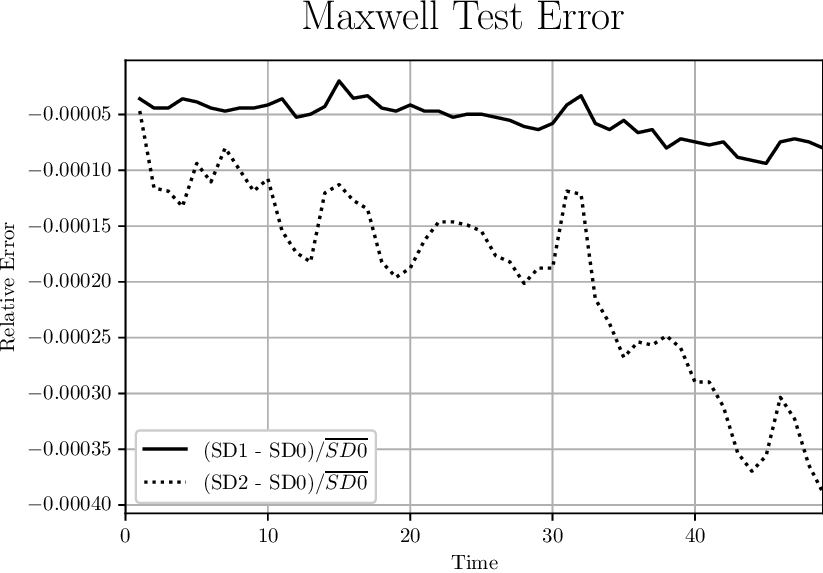}
\caption{Energy conservation error of the ELF propagation on the sphere. SD0: LOCAL moving frames with $\mathcal{G}=0$. Divergence conservation scheme (SD1) uses Eq. \eqref{Restorescheme1}. Connection preservation scheme (SD2) uses Eq.\eqref{Restorescheme2}. }
\label {maxwellerror}
\end{figure}

For the test problem, consider the initial shock of the electromagnetic pulse and energy conservation error during the extremely low frequency (ELF) propagation on the sphere. The TM mode is chosen. A shock is initiated at $(1/\sqrt{2},0,1/\sqrt{2})$ to propagates infinitely between the initial point and antipode without losing energy. Fig. \ref{maxwellerror} presents the difference in the energy conservation error for the two schemes versus time. The SD2 scheme noticeably shows better conservation error than does the SD2 scheme. Similar to the static curl operator test, the covariant derivative component $\nabla \times \mathbf{e}^3$ is computed more accurately in the connection preservation scheme than the divergence preservation scheme for the deformed domain.

\section{Application to SWEs}
As an application of the proposed scheme, SWEs is chosen on the sphere as the governing equation is given as follows:
\begin{linenomath}
\begin{align}
& \frac{\partial H}{\partial t} + \nabla \cdot (H \mathbf{u} ) = 0 ,  \label{SWE1} \\
& \frac{\partial H \mathbf{u} }{\partial t} + \nabla \cdot (H \mathbf{u} \mathbf{u}  ) + \frac{g}{2} \nabla H^2 = f H ( \mathbf{u} \times \mathbf{k} ) + g H \nabla H_0, \label{SWE2}
\end{align}
\end{linenomath}
where $\eta$ and $\mathbf{u}$ represent the free surface elevation and the depth-averaged velocity, respectively. $H_0$ is the still water depth and $H= H_0 + \eta$ is the total water depth. $f$ and $g$ are the Coriolis parameter and the gravitational constant, respectively. The divergence preservation scheme is obtained similarly as for other PDEs. However, the connection preservation scheme requires further modifications in the SWEs because of the (1) zero tendency condition, (2) gradient computation, and (3) differentiation of different vectors.

\subsection{Zero tendency condition}
Let $\mathbf{e}^i$ be moving frames that are constructed at each point on the sphere. Consider $\mathbf{u} = u \mathbf{e}^1 + v \mathbf{e}^2$. Multiply $\mathbf{e}^i$ with Eq. \eqref{SWE2} to obtain the following expression: 
\begin{linenomath}
\begin{align}
& \frac{\partial H \mathbf{u} }{\partial t} \cdot \mathbf{e}^1 + \nabla \cdot (H u \mathbf{u}  ) + H \mathbf{u} \cdot \nabla \mathbf{e}^1 \cdot \mathbf{u}+ \frac{g}{2} \nabla H^2 \cdot \mathbf{e}^1 = f H u + g H \nabla H_0 \cdot   \mathbf{e}^1, \label{SWE21} \\
& \frac{\partial H \mathbf{u} }{\partial t} \cdot \mathbf{e}^2 + \nabla \cdot (H v \mathbf{u}  ) + H \mathbf{u} \cdot \nabla \mathbf{e}^2 \cdot \mathbf{u} + \frac{g}{2} \nabla H^2 \cdot \mathbf{e}^2 = - f H v + g H \nabla H_0 \cdot  \mathbf{e}^2 . \label{SWE22} 
\end{align}
\end{linenomath}
The above equations are obtained by assuming the following zero tendency condition:
\begin{linenomath}
\begin{equation}
\mathbf{u} \cdot \nabla \mathbf{e}^i \cdot \mathbf{u} = 0, \hspace{1cm} 1 \le i \le 2. \label{zerotend} 
\end{equation}
\end{linenomath}
If moving frames are aligned in a more \textit{Euclidean connection}, the tendency is closer to zero. In other words, if the moving frames' curvature within an element is closer to zero, the left hand side in Eq.\eqref{zerotend} is also closer to zero. For the Williamson's tests, the zero tendency of the spherical coordinate is assumed to be satisfied, even when the left hand side in Eq.\eqref{zerotend} is not necessarily zero along the spherical coordinate axis. the left hand side in Eq.\eqref{zerotend} is implicitly assumed to be zero or at least approximately negligible compared with discretization error. Because any Euclidean axis or spherical coordinate axis provides zero tendency conditions, only LOCAL and LOCSPH satisfy the zero-tendency condition.

\subsection{Gradient computation}
Another problem pertains to formulating the gradient in the SWEs for the divergence/connection conservation scheme. In the SWEs, two gradients emerge: $ (g/2) \nabla H^2$ and $g H \nabla H_0$ in Eq.\eqref{SWE2}. In most numerical schemes \cite{Legat,Giraldo2005,Lauter2008,Nair2005}, the gradient of $H^2$, or $\nabla H^2$, is computed in the weak formulation. For the given moving frames $\mathbf{e}^i$, the weak formulation of $\nabla H^2$ can be computed similarly to the divergence as follows.

\begin{align}
\int \nabla H^2 \cdot \mathbf{e}^i \varphi dx = - \int \nabla \varphi \cdot \mathbf{e}^i H^2 dx - \int H^2 (\nabla \cdot \mathbf{e}^i ) \varphi d x + \int_{\partial \Omega} H^2 (\mathbf{e}^i \cdot \mathbf{n} ) \varphi d s  .  \label{dirderiv}
\end{align}
Note that spurious divergence is not added to the weak formulation of the gradient, even though the integration by parts of its weak formulation is similar to that of divergence. We assume that no spurious divergence is generated for the gradient in the slightly deformed domain. Several computational tests have confirmed this claim but are not included in the paper.

Another gradient, or $g H \nabla H_0$, is computed directly because the integration of part for $H \nabla H_0$ is impossible when $H_0$ is not a constant. We observed that the existence of this direct computation of the gradient causes deteriorated accuracy for moving frames other than LOCAL moving frames. 

\subsection{Differentiation of vectors}
In the SWEs in Eqs. \eqref{SWE1} and \eqref{SWE2}, two different differentiations of a vector occur: the first occurs for divergence computations such as
\begin{equation}
\nabla \cdot (H \mathbf{u}), ~~~\nabla \cdot (H u \mathbf{u}), ~~~\nabla \cdot (H v \mathbf{u}) \label{Diffvec1}
\end{equation}
The second occurs for rotation effect to yield the following component non-trivial in Eqs. \eqref{SWE21} - \eqref{SWE22}.
\begin{equation}
H \left (   u \frac{\partial \mathbf{e}^1 }{\partial t} +  v \frac{\partial \mathbf{e}^2 }{\partial t}  \right ) \cdot \mathbf{e}^i,~~~~~~~~ 1 \le i \le 2  \label{Diffvec2}
\end{equation}
The derivative $\partial \mathbf{e}^i / \partial t$ is another direct differentiation with respect to moving frames $\mathbf{e}^i$ such as
\begin{equation*}
\frac{\partial \mathbf{e}^i}{\partial t} = \sum_{k=1}^3 \left [ u ( \nabla e_k^i \cdot \mathbf{e}^1 ) + v ( \nabla e_k^i \cdot \mathbf{e}^2 ) \right ] \mathbf{x}^k,
\end{equation*}
where $\mathbf{x}^k$ is the $k$th Cartesian frame. For the computation of Eqs. \eqref{Diffvec1} and \eqref{Diffvec2}, the choice of moving frames should be theoretically independent. However, this is not true for a slightly deformed domain $\Omega^{\varepsilon}$. Adopting LOCAL moving frames for $\Omega^{\varepsilon}$ reflects the current geometry, or connection. For the connection preservation scheme, LOCSPH moving frames, which are the original connection of $\Omega^0$, are adopted for the expansion of the vector in Eqs. \eqref{Diffvec1} and \eqref{Diffvec2}.

\subsection{Divergence/connection preservation scheme}
Considering the aforementioned three factors for the improved accuracy of the numerical scheme, the following weak formulation with moving frames is proposed for SWE equations in the context of the discontinuous Galerkin method. 
 \begin{align}
&\int \frac{\partial H}{\partial t} \varphi dx  + \Delta_1 + \Delta_2   - \int \mathcal{G} (H\mathbf{u})  \varphi dx    = 0,  \label{SWEint1}  \\ 
 & \int\frac{\partial H u  }{\partial t} \varphi dx + \Gamma_1  + \Gamma_2  + \Gamma_3  + \Gamma_4   + \Gamma_5  -  \int \mathcal{G}(H u \mathbf{u}) \varphi dx    \nonumber \\
 &~~~~~~~~~~~~~~~~~~~~~~~~~~~~~~~~~~~~~~~~~~~~~=  \int_{\Omega}  f \left( H v \right) \varphi dx + \Gamma_6   . \label{SWEint2}   \\
    & \int\frac{\partial H v  }{\partial t} \varphi dx + \Lambda_1 + \Lambda_2  + \Lambda_3  + \Lambda_4 + \Lambda_5   -   \int \mathcal{G}(H v \mathbf{u}) \varphi dx     \nonumber \\
  &~~~~~~~~~~~~~~~~~~~~~~~~~~~~~~~~~~~~~~~~~~~~~   =  - \int_{\Omega}  f \left( H u \right) \varphi dx + \Lambda_6 \label{SWEint3}  .
\end{align}
where
\begin{align*}
\Delta_1 &=  - \int  H  u^d \mathbf{d}^1  \cdot  \nabla \varphi  dx - \int  H  v^d \mathbf{d}^2  \cdot  \nabla \varphi  dx, \\
\Delta_2 &=  \int_{\partial \Omega} \widetilde{H} \tilde{\mathbf{u}} \cdot {\mathbf{n}} \varphi ds,
\end{align*}
and
\begin{align*}
\Gamma_1 &=  \int   H  \left ( u^d  \frac{\partial \mathbf{d}^1 }{\partial t} + v^d  \frac{\partial \mathbf{d}^2 }{\partial t}  \right ) \cdot \mathbf{e}^1 \varphi \,dx  \\
\Gamma_2 &=  -  \int   H {u} (u^d \nabla \varphi \cdot \mathbf{d}^1  + v^d \nabla \varphi \cdot \mathbf{d}^2  )  d x, \\
\Gamma_3 &=  -  \int \frac{g H^2}{2}   \nabla \varphi \cdot \mathbf{e}^1 d x, \\
\Gamma_4 &=   \int_{\partial \Omega}  \left [ \widetilde{H}  \widetilde{u}^2  + \frac{g}{2} \widetilde{H}^2  \right ] (\mathbf{e}^1 \cdot \mathbf{n}) \varphi dx, \\
\Gamma_5 &= - \int   \frac{g H^2}{2 }  \left ( \nabla \cdot \mathbf{e}^1 \right )  \varphi \,dx, \\
\Gamma_6 &=  \int  g H \nabla H_0 \cdot  \mathbf{e}^1 \varphi\, dx,
\end{align*}
and
 \begin{align*}
\Lambda_1 &= \int   H  \left ( u^d  \frac{\partial \mathbf{d}^1 }{\partial t} + v^d  \frac{\partial \mathbf{d}^2 }{\partial t}  \right ) \cdot \mathbf{e}^2 \varphi \,dx, \\
\Lambda_2 &= -  \int   H {v} (u^d \nabla \varphi \cdot \mathbf{d}^1  + v^d  \nabla \varphi \cdot \mathbf{d}^2  )  d x , \\
\Lambda_3 &=  -  \int \frac{g H^2}{2}   \nabla \varphi \cdot \mathbf{e}^2 d x, \\
\Lambda_4 &=  \int_{\partial \Omega}  \left [ \widetilde{H}  \widetilde{u}^2  + \frac{g}{2} \widetilde{H}^2  \right ] (\mathbf{e}^2 \cdot \mathbf{n}) \varphi dx , \\
\Lambda_5 &=  - \int   \frac{g H^2}{2 }  \left ( \nabla \cdot \mathbf{e}^2 \right )  \varphi \,dx , \\
\Lambda_6 &=  \int  g H \nabla H_0 \cdot  \mathbf{e}^2 \varphi\, dx .
 \end{align*}
 The moving frames $\mathbf{d}^i$ are also orthonormal vectors. However, they are not necessarily the same as $\mathbf{e}^i$. The moving frames $\mathbf{d}^i$ preserve the original connection of the original domain $\Omega^0$. In our scheme, $\mathbf{e}^i$ is chosen as LOCAL moving frames and $\mathbf{e}^i_c$ is chosen as LOCSPH moving frames. For example we have the following equality in any domain. 
 \begin{equation*}
 \mathbf{u} =  u \mathbf{e}^1 + v \mathbf{e}^2 = u^d \mathbf{d}^1  + v^d \mathbf{d}^2  .
 \end{equation*}

Without the spurious divergence, the numerical scheme in Eqs. \eqref{SWEint1} - \eqref{SWEint3} is exactly the same as the typical numerical scheme with moving frames, as used in ref. \cite{MMF3}. Note that the $\mathcal{G}$ component has a different velocity vector as indicated in the parentheses. If LOCAL is used for $\mathbf{d}^i $, then the numerical scheme is a divergence preservation scheme, the same as in Eq. \eqref{Restorescheme1}. If LOCSPH is used for $\mathbf{d}^i $, then the numerical scheme is both a divergence and a  connection preservation scheme. We call this scheme as \textit{divergence/connection preservation scheme}, as described in detail as follows.

Consider a velocity vector that lies on the slightly deformed mesh $\Omega^{\varepsilon}$, but is expanded by LOCSPH, denoted by $  {\mathbf{v}}^{0}_{\| \varepsilon}$. The subscript $\| \varepsilon$ implies that the corresponding quantity is projected onto $\Omega^{\varepsilon}$. Similarly, the del operator is expanded by LOCSPH. This means that the differentiation occurs in the direction of the LOCSPH frame and is projected at the same time because the grid points lie on the slightly deformed domain $\Omega^{\varepsilon}$. By substituting $\mathbf{F} = \mathbf{k}^{\varepsilon} \times  {\mathbf{v}}^{0}_{\| \varepsilon}$ into the following Stokes theorem, we obtain

\begin{equation}
\iint_{\Omega^{\varepsilon}_e} ( \nabla^{0}  \cdot \mathbf{v}^{0} )  d A = \iint_{\Omega^{0 \varepsilon}_e} ( \nabla^{0}_{\| \varepsilon}  \cdot \mathbf{v}^{0}_{\| \varepsilon} )  d A -  \iint_{\Omega^{\varepsilon}_e }  \mathcal{G}(\mathbf{k}^{\varepsilon},\nabla^{0}_{\| \varepsilon},\mathbf{v}^{0}_{\| \varepsilon})  d A  . \label{Restorescheme3}
\end{equation}
Eq. \eqref{Restorescheme3} is fundamentally a divergence preservation scheme because it retrieves the original divergence $\nabla^0 \cdot \mathbf{v}^0$. At the same time, the use of the original connection $\nabla^0$ preserves the original differentiation connection.

For the numerical tests in the following subsections, the Lax-Friedrich flux is used for the numerical flux. The classical explicit fourth-order Runge-Kutta scheme with a sufficiently small $\Delta t$ is used for time marching. The details of the numerical flux and scheme are well described in ref. \cite{MMF4}. The numeric test problems are also described in ref. \cite{MMF4,Williamson1992}, thus they are not presented in this paper.

\begin{figure}[ht]
\begin{center}
\resizebox{1.0\textwidth}{!}{%
\subfloat[Initial condition] {\label{SWESteadyInit}\includegraphics[height=4cm]{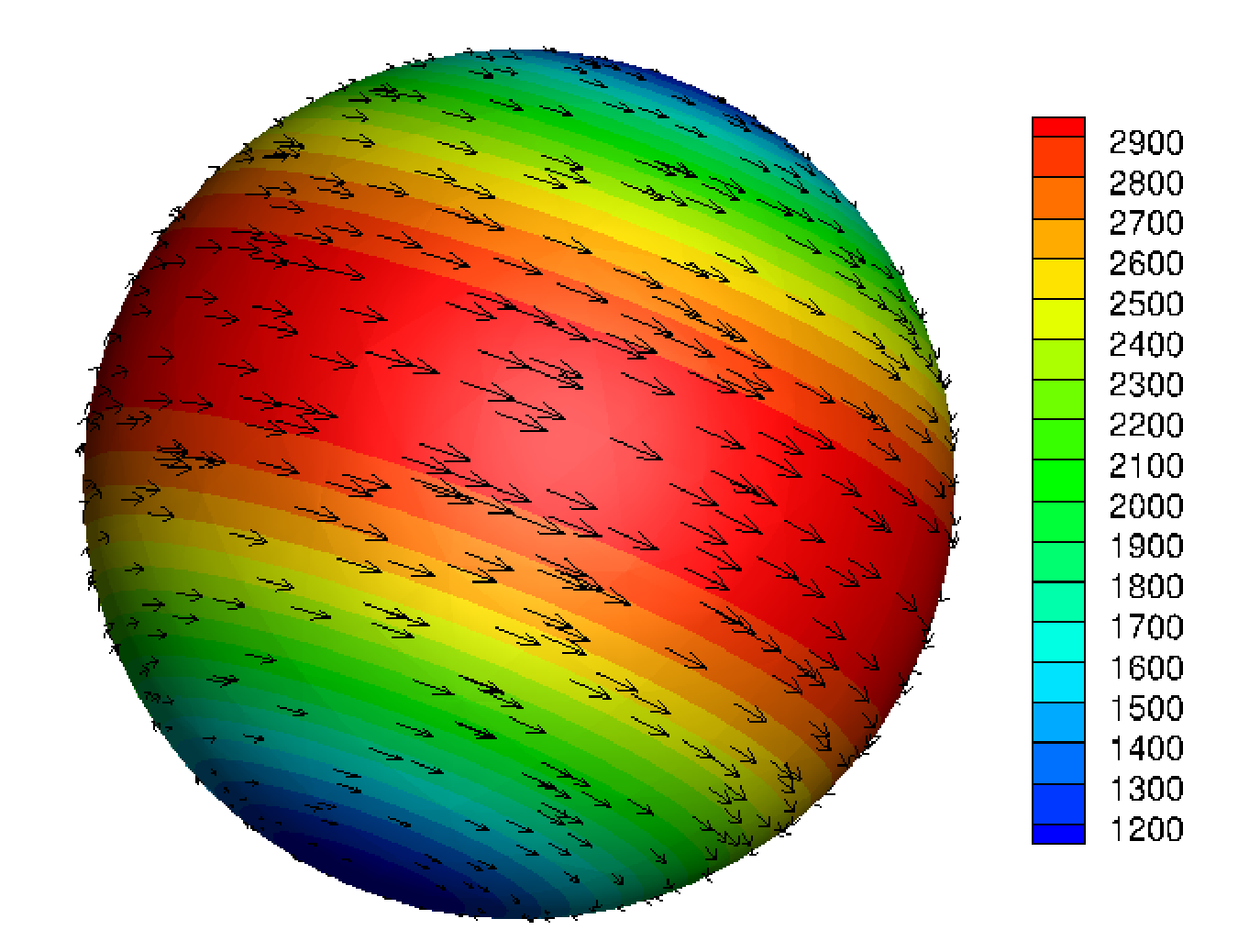} }
\quad
\subfloat[$L_2$ error] {\label{SWESteadyL2}\includegraphics[height=5cm]{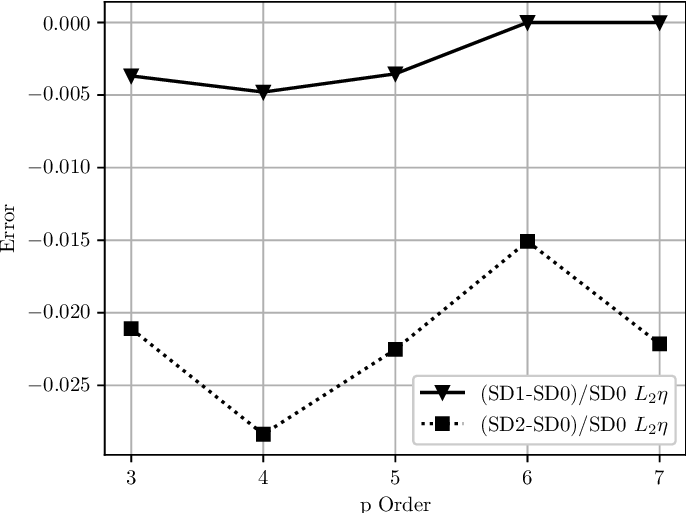} } 
\quad
\subfloat[Mass conservation error] {\label{SWESteadyMass}\includegraphics[height=5cm]{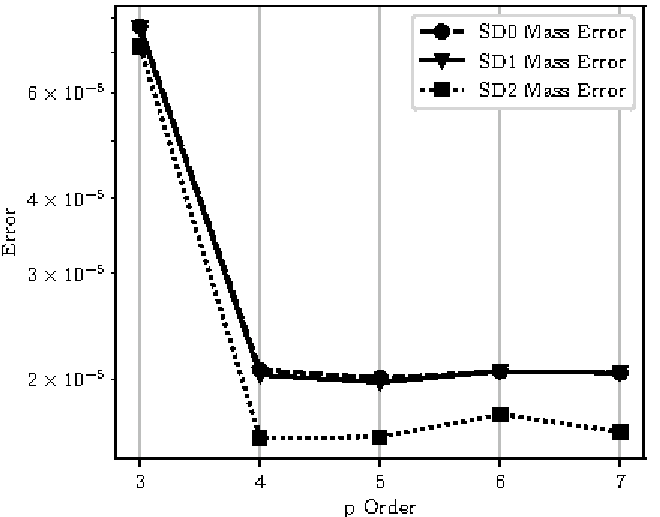} }
}
\end{center}
\caption{Steady zonal flow: Initial condition of $\eta$ and velocity vector, overall relative $L_2$ error, and mass conservation error for $3 \le p \le 7$. $h=0.4$ with 480 elements. Computed up to $T$=$5$ (days). SD0: LOCAL moving frames with $\mathcal{G}=0$. Divergence conservation scheme (SD1) uses Eq. \eqref{Restorescheme1}. Divergence/connection preservation scheme (SD2) uses Eq.\eqref{Restorescheme3}. }
\label {SWESteadyZonal}
\end{figure}

\subsection{Steady zonal flow}

The first test is the steady zonal test with constant free surface elevation and velocity vector over time. The angle of the velocity vector with respect to the polar axis is set to $\pi/4$. The still water depth ($H_0$) is chosen as $1.0$ and $\nabla H_0 = 0$. Error is measured at $T=5.0$ (days).

Figure. \ref{SWESteadyZonal} displays the overall relative $L_2$/$L_{\infty}$ error and the mass/energy conservation error compared with the initial mass and energy. The overall errors are saturated for $p \ge 4$ in all the schemes. However, the overall $L_2$ error for the SD2 scheme is significantly smaller than those of SD0 and SD1. Fig. \ref{SWESteadyL2} shows that the relative error of SD2 in comparison with SD1 is at least 50$\%$ smaller than that of SD1. The overall error increases as $T$ increases.

Similarly, mass and energy conservation error are both saturated for $p \ge 4$ in all the schemes. Conservation error is more severe than the overall error. The divergence/connection preservation scheme (SD2) shows noticeably smaller errors in comparison with SD0 and SD1. For mass conservation error, the difference between the SD2 scheme and SD0/SD1 is approximately $3.00 \times 10^{-7}$ at $T=5$ for $p\ge4$. The error difference is small but significant when the mesh error is approximately in the range of $2.00 \times 10^{-5} \sim 3.00 \times 10^{-4}$. However, this small difference contributes to a meaningful deterioration of the overall error in long time integration.

\subsection{Unsteady zonal flow}
Contrary to the steady zonal test, the unsteady zonal test in the Williamson's test suites uses a variable still surface elevation of $H_0$ as follows.
\begin{equation*}
H_0 (\theta) = \frac{133681}{r_a g} - \frac{10}{r_a g} - \frac{(\Omega_f \sin \theta)^2}{2g} ,
\end{equation*}
where $r_a$ is the radius of the earth, $g$ is the gravitational constant, and $\Omega_f$ is the Earth's rotational frequency. Note that $g H \nabla H_0 \cdot \mathbf{e}^i$ should be computed directly because it cannot be computed in a weak form by integration by part. Thus, the divergence/connection scheme cannot be applied to $H \nabla H_0$.

\begin{figure}[ht]
\begin{center}
\resizebox{1.0\textwidth}{!}{%
\subfloat[Initial condition] {\label{SWEUnsteadyInit}\includegraphics[height=4cm]{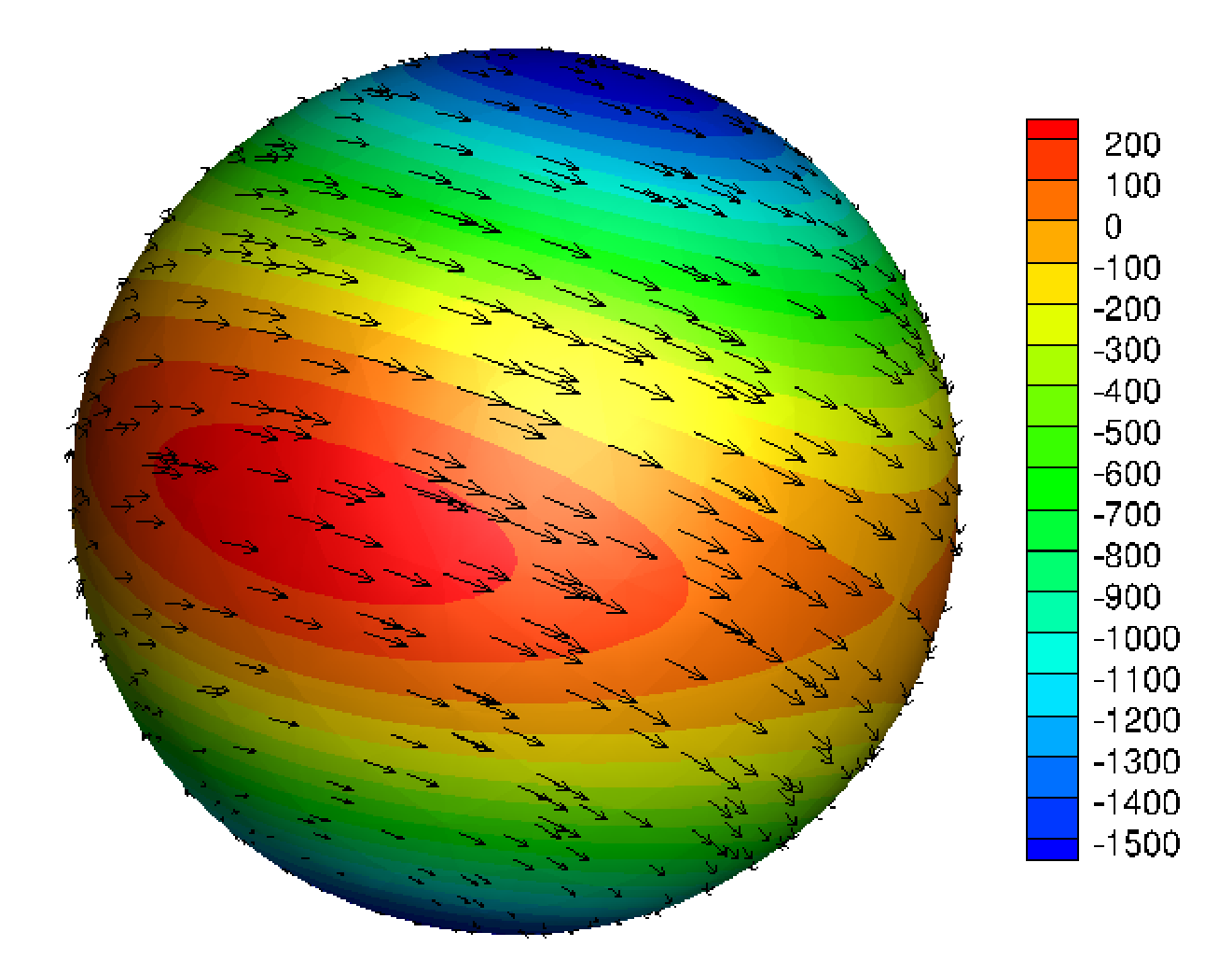} }
\quad
\subfloat[$L_2$ error] {\label{SWEUnsteadyL2}\includegraphics[height=5cm]{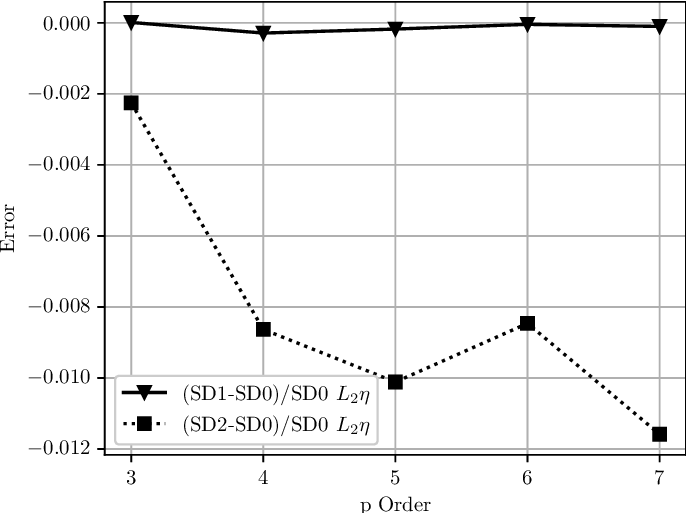} } 
\quad
\subfloat[Mass conservation error] {\label{SWEUnsteadyMass}\includegraphics[height=5cm]{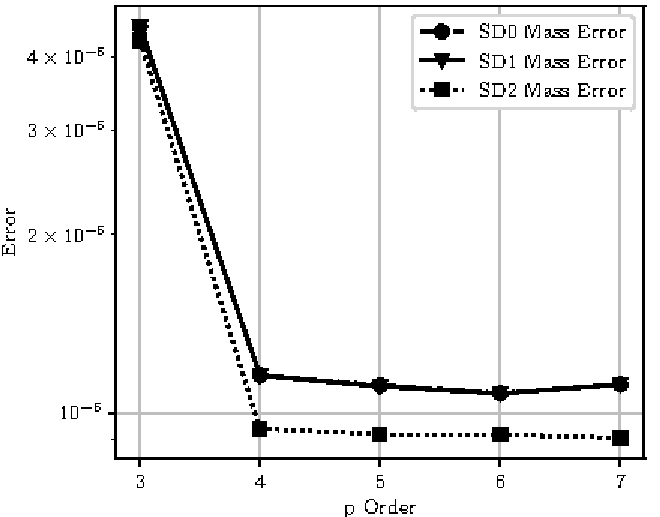} }
}
\end{center}
\caption{Unsteady zonal flow: Initial condition of $\eta$ and velocity vector, overall relative $L_2$ error, and mass conservation error for $3 \le p \le 7$. $h=0.4$ with 480 elements. Computed up to $T$=$5$ (days). SD0: LOCAL moving frames with $\mathcal{G}=0$. Divergence conservation scheme (SD1) uses Eq. \eqref{Restorescheme1}. Divergence/connection preservation scheme (SD2) uses Eq.\eqref{Restorescheme3}. }
\label {SWEUnsteadyZonal}
\end{figure}

Similar to the error plot of the steady zonal test in Fig. \ref{SWESteadyZonal}, Fig. \ref{SWEUnsteadyZonal} presents the overall error and mass conservation error for the SD0, SD1, and SD2 schemes. The overall errors are also saturated for $p \ge 4$ in all the schemes. The overall $L_2$ error for the SD2 scheme is slightly better than those of SD0 and SD1 because of the relatively short time integration up to $0.5$(days). The mass conservation error difference is approximately $3.30 \times 10^{-7}$ at $T=5$ for $p\ge4$, similar to the result of the steady zonal flow test. Adopting LOCAL moving frames for $\mathbf{e}^i$ seems to be the primary reason for the accuracy maintained for the direct derivative of the non-constant $H_0 (\theta)$, regardless of $\mathbf{d}^i$.

\subsection{Rossby--Haurwitz flow}
\begin{figure}[ht]
\begin{center}
\resizebox{1.0\textwidth}{!}{%
\subfloat[Initial condition] {\label{SWERossbyInit}\includegraphics[height=4cm]{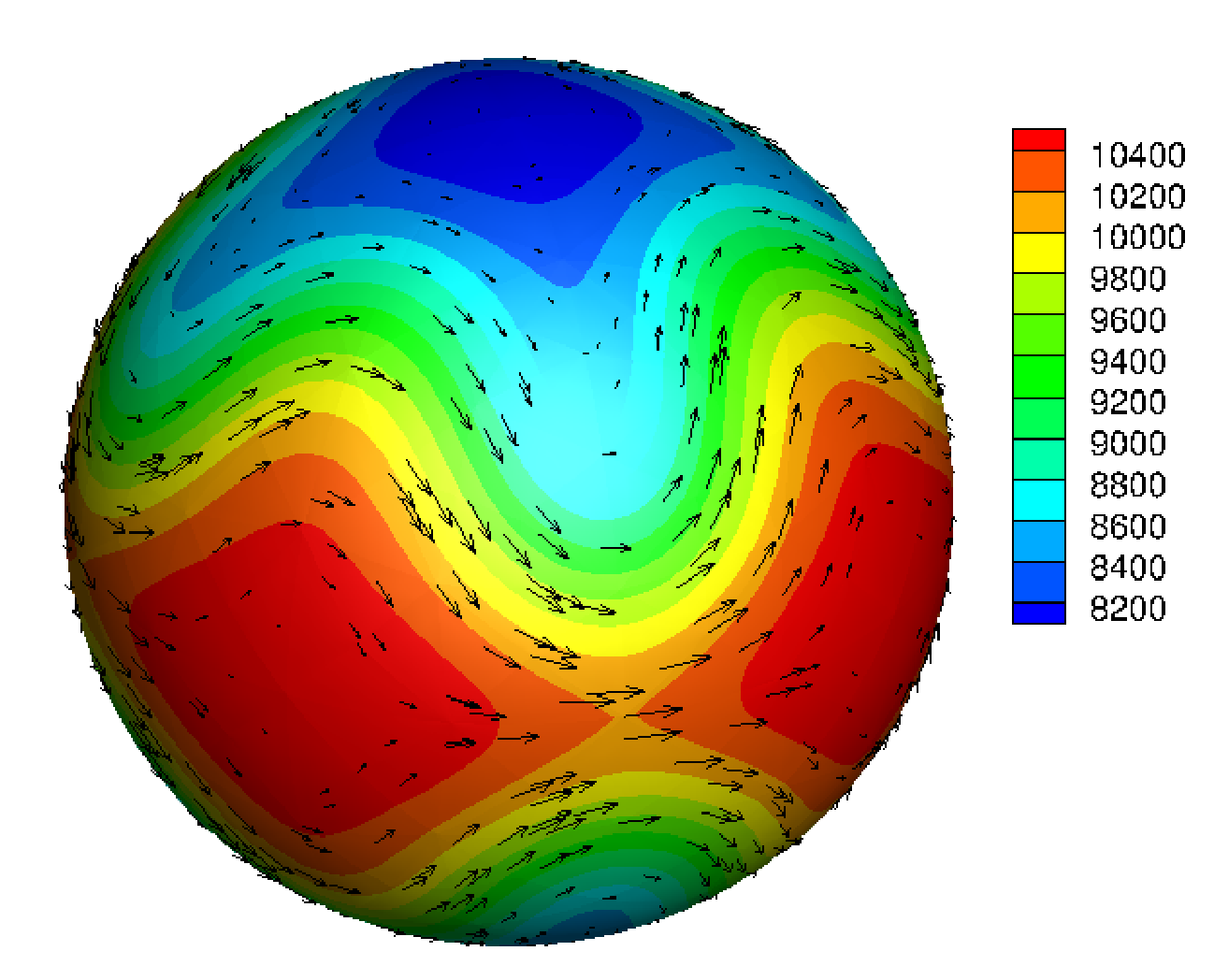} }
\quad
\subfloat[Mass conservation error] {\label{SWEUnsteadyMass}\includegraphics[height=5cm]{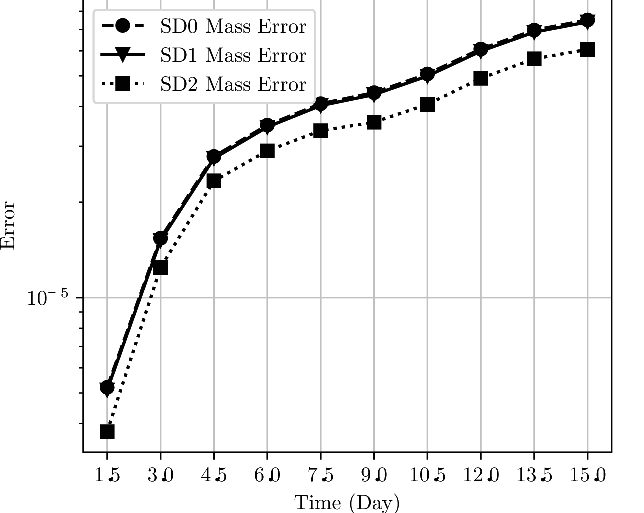} }
\quad
\subfloat[Energy conservation error] {\label{SWEUnsteadyMass}\includegraphics[height=5cm]{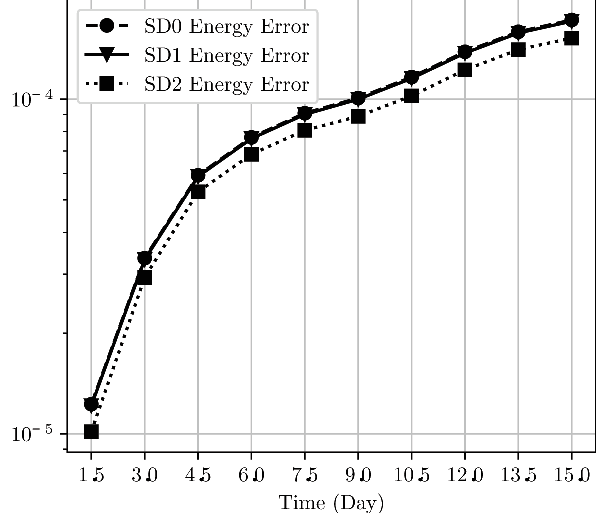} }
}
\end{center}
\caption{Rossby-Haurwitz flow simulation: Initial condition of $\eta$ and velocity vector, mass, and energy conservation error are displayed. $p=5$. $h=0.4$ with 480 elements. Computed up to $T=15.0$ (days). SD0: LOCAL moving frames with $\mathcal{G}=0$. Divergence conservation scheme (SD1): Use Eq. \eqref{Restorescheme1}, Connection preservation scheme (SD2): Use Eq.\eqref{Restorescheme3}. }
\label {SWERossby}
\end{figure}

The third test problem is the Rossby--Haurwitz flow on the sphere. There is no exact solution to this problem. However, the periodic distribution of the free surface elevation ($\eta$) and the velocity vector ($\mathbf{u}$) along the $\phi$-axis rotate around the axis of the sphere. The conservation error of mass and energy is considered the most important measure for demonstrating the accuracy of the implemented numerical schemes. The initial distribution of $\eta$ and $\mathbf{u}$ is more complicated than those of the steady zonal flow, as shown in Fig. \ref{SWERossby}. Let $H_0$ be constant, set to $1.0$, and $\nabla H_0$ be zero. Thus, the LOCSPH with the $\mathcal{G}$ scheme can be implemented for better accuracy. The SWEs equations are computed up to $T=15$ days.

Fig. \ref{SWERossby} demonstrates the mass and energy conservation error versus time. The SD2 scheme clearly shows better accuracy both in mass and energy conservation error than do the SD0 and SD1 schemes. The SD1 scheme's mass and energy conservation errors are more accurate than those of SD0 by $1.09\times 10^{-06}$ and $1.51\times 10^{-06}$, respectively. The SD2 scheme's mass and energy conservation errors are more accurate than those of SD0 by $1.45\times 10^{-05}$ and $2.02\times 10^{-05}$, respectively. The SD2's improved accuracy is significant in comparison with the mesh error in the range of $2.00 \times 10^{-5} \sim 3.00 \times 10^{-4}$.

\subsection{Steady jet flow}

\begin{figure}[ht]
\begin{center}
\resizebox{1.0\textwidth}{!}{%
\subfloat[Initial condition] {\label{SWEjetInit}\includegraphics[ width=4cm]{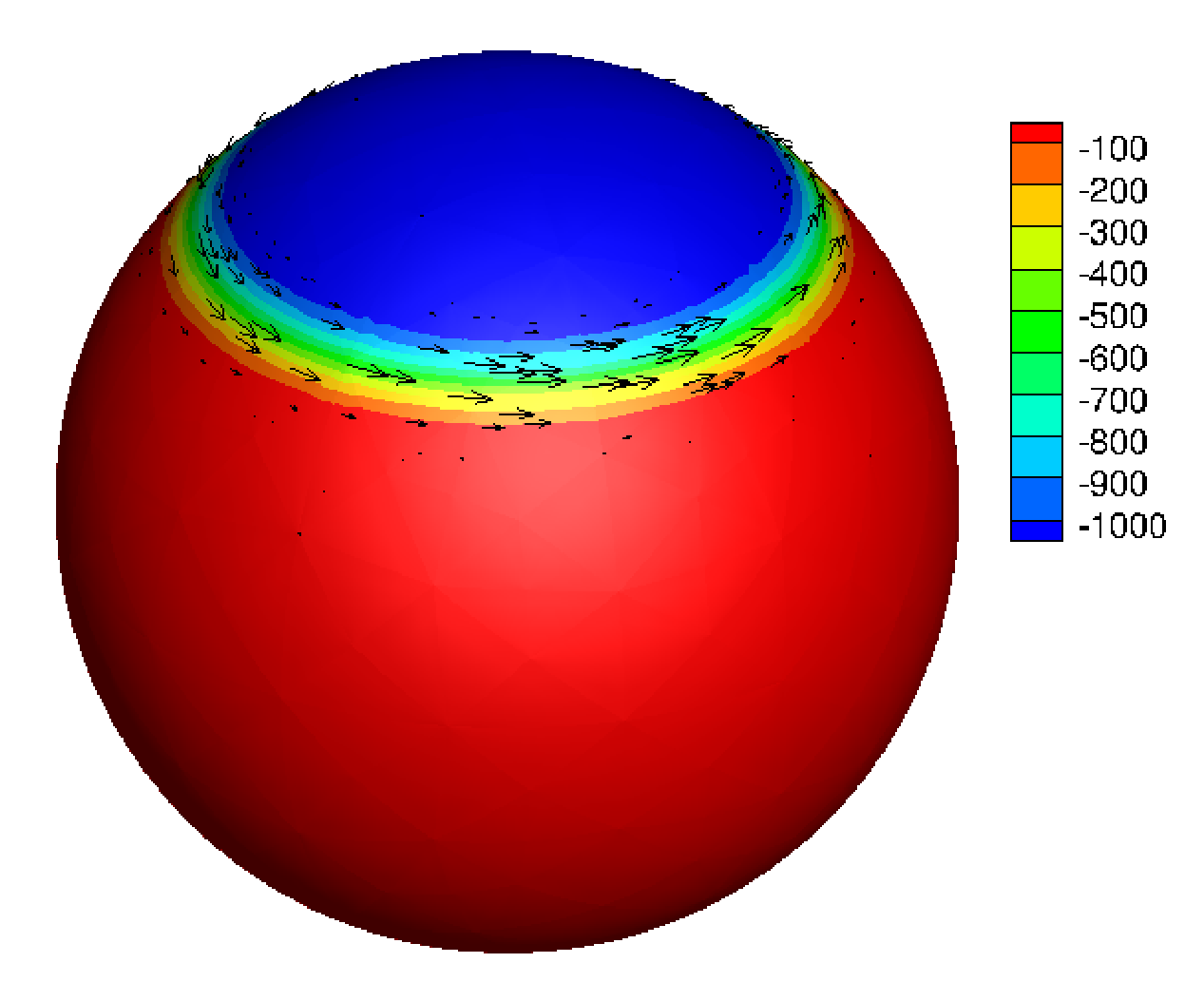} }
\subfloat[Mass conservation error] {\label{SWEjetError}\includegraphics[width=5cm]{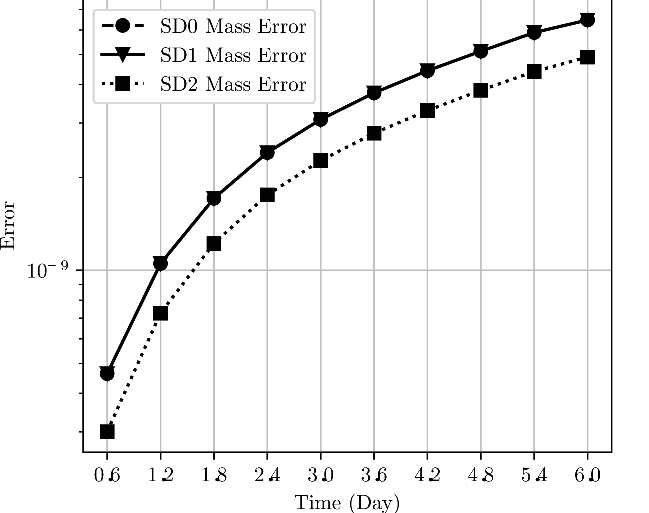} }}
\end{center}
\caption{Jet flow: Initial condition and mass conservation error. $p=5$. $h=0.08$ with 10752 curved elements.}
\label {SWEJet}
\end{figure}

As the last test problem with the constant $H_0$, the jet flow simulation is known to be one of the most challenging tests as the initial condition is presented in Fig. \ref{SWEjetInit}. One difficulty is that a very fine spatial resolution should be provided for accurate computation. For the weak formulation, it is known that the edge of curved elements should be less than $h$=$0.1$. Consequently, the number of curved triangular elements for a sphere should be at least more than $10,000$ for most numerical schemes. For our simulation, 10752 curved elements are used with $p=5$, equivalent to approximately $4.5\times 10^5$ grid points. Perturbed jet simulations require more space resolutions. Fig. \ref{SWEjetError} demonstrates the improved accuracy in mass conservation in a steady jet flow simulation up to $T=6$ (days). At $T=6$, the mass conservation error of the SD1 scheme is $99.97\%.$ compared to that of the SD0 scheme. However, the mass conservation error of the SD2 scheme is only $75.87\%$ compared to that of the SD0 scheme. Results are similar for perturbed jet simulations for increased spatial resolution for $p>5$. Thus, this result is omitted in the paper.

\section{Discussion}
This paper proposes new schemes to improve the overall accuracy and particularly the conservation error, in a slightly deformed curved surface. The improvements in accuracy are made possible by retrieving the original divergence or the original connection from the original domain, by using the original surface normal vector or the original del operator or connection. When the discretization error is sufficiently low, the proposed schemes provide significantly improved accuracy for the overall error and mass/energy conservation. The main drawback of the proposed schemes is that extra computational costs are required to compute the spurious divergence $\mathcal{G}$ at every time step. Accuracy gain is relatively small considering the added computational costs. However, the slightly improved accuracy identifies the source of the conservation error by geometric approximation error. Thus, this scheme can significantly improve the overall error in long time integration where conservation is critical, i.e., climate simulations.

Moreover, the results of this study suggest that when the surface domain has been modified from the original domain, retrieving the changes in divergence could be less critical than retrieving the connection of the original domain. In other words, the source of errors for a slightly modified domain is more affected by the changes in the path for the differentiation of vectors or the so-called \textit{covariant derivative} of vectors. Retrieving the original connection requires prior knowledge of the domain, such as the surface normal vector $\mathbf{k}^0$ or the $\nabla$ operator $\nabla^0$. This result can be applied to various applications, for example, in divergence computation in data science, where the original information is slightly modified, but a theoretical configuration exists.

\section*{Acknowledgements}
This research was supported by the National Research Foundation of Korea (NRF-2021R1A2C109297811). The second author is supported by the R\&D project on ``Development of the Next-generation Operational System of the Korea Institute of Atmospheric Prediction Systems (KIAPS)'', funded by Korea Meteorological Administration (KMA2020-02213).

\bibliographystyle{amsplain}
\bibliography{DCPScheme_arXiv_2022}

\end{document}